\documentclass[10pt]{article}

\RequirePackage{amsmath}
\RequirePackage{amsfonts}
\RequirePackage{amssymb}
\RequirePackage{hyperref}
\RequirePackage{multirow}
\RequirePackage{algorithm}
\RequirePackage{algpseudocode}
\RequirePackage{graphicx}
\RequirePackage{xcolor}
\RequirePackage{amsthm}
\RequirePackage{amsmath}
\RequirePackage{amsfonts}
\RequirePackage{amssymb}
\RequirePackage{hyperref}
\RequirePackage{multirow}
\RequirePackage{algorithm}
\RequirePackage{algpseudocode}
\RequirePackage{graphicx}
\RequirePackage{xcolor}
\allowdisplaybreaks
\RequirePackage{enumitem}

\usepackage{dpr_fec,ifthen,dsfont}

\newtheorem{definition}{Definition}[section]
\newtheorem{lemma}{Lemma}[section]
\newtheorem{theorem}{Theorem}[section]

\usepackage{isetechreport}
\def\reportyear{26T}
\def\reportno{013}
\def\revisionno{0}
\def\originaldate{\today}
\def\revisiondate{\today}
\coraltrue
\cvcrfalse

\begin{document}

\title{A Local-Linearly Convergent Algorithm for Nonconvex Equality-Constrained Optimization}

\author{Frank E.~Curtis\thanks{E-mail: frank.e.curtis@lehigh.edu}}
\author{Lingjun Guo\thanks{E-mail: lig423@lehigh.edu}}
\author{Daniel P.~Robinson\thanks{E-mail: daniel.p.robinson@lehigh.edu}}
\affil{Department of Industrial and Systems Engineering, Lehigh University}
\titlepage

\maketitle

\begin{abstract}
  For solving nonconvex equality-constrained optimization problems, a recent Gradient-Eigenstep Algorithm by Goyens et al.~is an iteration-efficient approach, based on minimizing Fletcher's augmented Lagrangian function, for finding an approximate second-order stationary point from an arbitrary starting point. In this paper, the analysis of this algorithm is extended, offering a two-fold contribution. First, it is shown that a local-linear rate of convergence can be obtained by this method if it is initiated sufficiently close to a strong second-order stationary point and employs a sufficiently small step-size parameter and sufficiently large penalty parameter. In this case, the algorithm reduces to a gradient descent algorithm applied to minimize Fletcher's augmented Lagrangian. Second, as a particularly useful application of the first result, it is shown that the Gradient-Eigenstep algorithm can be used as an iteration-efficient subproblem solver in the context of a progressive sampling strategy for solving equality-constrained optimization problems when the objective and constraint functions are defined by large sample averages, ultimately offering an algorithm with an improved worst-case sample complexity when compared to an approach that solves a full-sample problem directly.
\end{abstract}

\section{Introduction}

This work studies the theoretical convergence-rate behavior of an algorithm for solving nonlinear equality-constrained continuous optimization problems. Problems of this type, and algorithms for solving them, have been the subject of research for decades. Many of the most successful algorithms can be characterized as being penalty methods or primal-dual methods; see, e.g., \cite[Chapter 17--18]{nocedal2006numerical}. These are local-search algorithms that seek a point satisfying either first- or second-order stationary conditions for optimality.

Our focus in this paper is on a particular penalty method.  Penalty methods solve a constrained problem by introducing a measure of constraint violation and adding it to the objective function, weighted by a penalty parameter.  The aim is to solve the original constrained problem through unconstrained optimization techniques, where the penalty parameter is adjusted in a principled manner to ensure that the constraints are eventually satisfied, at least in the limit.  The penalty function used in a penalty method can be characterized as being either \textit{inexact} or \textit{exact}.  An inexact penalty function is one for which the penalty parameter might need to be driven to infinity to ensure that constraint satisfaction is achieved in an asymptotic sense. By contrast, an exact penalty function is one for which, for a given problem, there exists a threshold for the penalty parameter, above which there is a correspondence between stationary points for the original constrained problem and stationary points for the penalty function. This feature of exact penalty functions is attractive, but an issue with many exact penalty functions is that they are nonsmooth.

Fletcher~\cite{fletcher1970class} provides an exact penalty function---which came to be known as Fletcher's augmented Lagrangian function---that is both exact and smooth. After it was introduced and shown to be exact, several articles \cite{DiPillo1994,DiPilloGrippo1986,DiPilloGrippo1989} explored the potential use of it in practical penalty methods for solving constrained problems, focusing in particular on its computational costs related to the computation of least-squares Lagrange multipliers and their derivatives. Recently, the article~\cite{estrin2020} proposed an algorithm based on Fletcher's augmented Lagrangian that only requires a linear system solve in each iteration. Even more recently, \cite{GoyeEfteBoum2024} extended the work in \cite{DiPilloGrippo1986} by proposing an algorithm that exploits potential negative curvature directions to obtain an approximate second-order stationary point with a complexity guarantee.

Our contributions in this paper are twofold. First, we show that if an initial point is given that is sufficiently close to a strong second-order stationary point satisfying a regularity condition, then if the penalty parameter is sufficiently large and the step sizes are chosen sufficiently small, then an application of gradient descent to minimize Fletcher's augmented Lagrangian function can be used to obtain a linear rate of decrease in the norm of Fletcher's augmented Lagrangian, ultimately yielding an even more accurate approximate strong second-order stationary point. In this pursuit, we show that the Gradient-Eigenstep method in \cite{GoyeEfteBoum2024} reduces to gradient descent with a local-linear rate of convergence in the vicinity of certain strong second-order stationary points. Second, as a particularly useful application of these first results, we prove strong worst-case complexity properties for a recently proposed progressive sampling method for solving certain equality-constrained problems; see \cite[Algorithm 1]{CurtGuoRobi2025}. In particular, the progressive samping method is designed for solving problems where the objective and constraint functions are formed by large-scale sample average approximations of functions defined by expectations. The idea of progressive sampling is to solve a subsampled problem to modest accuracy, and then use the solution as a starting point for solving a subsequently subsampled problem to higher accuracy, progressively increasing the sample size until the original sample average approximation problem is solved to a desired accuracy. If the original problem satisfies a so-called strong-Morse property and the subsampled problems are solved through an application of gradient descent on Fletcher's augmented Lagrangian function with appropriate parameter choices, we show that the progressive sampling strategy yields an improved worst-case sample complexity property compared to a method that solves the full sample average problem directly.

\subsection{Notation}\label{sec.notation}

We use $\R{}$ to denote the set of real numbers, $\R{}_{\geq r}$ (resp.,~$\R{}_{>r}$) to denote the set of real numbers greater than or equal to (resp.,~greater than) $r \in \R{}$, $\R{n}$ to denote the set of $n$-dimensional real vectors, and $\R{m \times n}$ to denote the set of $m$-by-$n$-dimensional real matrices.  We define $\N{} := \{0,1,2,\dots\}$, and, for any integer $N \geq 1$, we use $[N]$ to denote $\{1, \dots, N\}$. For any finite set $\Scal$, we use $|\Scal|$ to denote its cardinality. For vectors, we define $\|\cdot\| := \|\cdot\|_2$, or else specify the norm explicitly.  We use $\|\cdot\|$ to denote the spectral norm of any matrix.

For any $A \in \R{m \times n}$ we use $\sigma_i(A)$ for each $i \in [\min\{m,n\}]$ to denote its $i$th largest singular value, and for any symmetric $A \in \R{n \times n}$ we use $\lambda_i(A)$ for each $i \in [n]$ to denote its $i$th largest eigenvalue.  Given any such $A$, we use $\Null(A)$ to denote its null space, i.e., $\{d \in \R{n} : Ad = 0\}$.  Assuming that $B \in \R{n \times m}$ has full-column rank, its Moore-Penrose pseudoinverse is $B^\dag := (B^TB)^{-1}B^T$.  Note that in this case of $B$ having full-column rank, the pseudoinverse $B^\dag$ is a left inverse of $B$ in the sense that $B^\dag B = I$. Given $B \in \R{n \times m}$ with full column rank, we use $\Rcal(B) := BB^\dag$ and $\Ncal(B^T) = I - \Rcal(B)$ to denote orthogonal projection matrices onto the span of the columns of $B$ and the null space of $B^T$, respectively.  That such $\Rcal(B)$ and $\Ncal(B^T)$ are orthogonal projection matrices follows since each of them are both symmetric and idempotent.

\subsection{Organization}

In \S\ref{sec.specific}, we show that the Gradient-Eigenstep method (\cite[Algorithm 1]{GoyeEfteBoum2024}) can reduce to a gradient-descent method that converges linearly in the neighborhood of certain second-order stationary points of an equality constrained optimization problem.  We apply this result to prove strong worst-case complexity properties of a progressive sampling method for equality constrained optimization in~\S\ref{sec.progressive}.  A summary and concluding remarks are given in~\S\ref{sec.conclusion}.

\section{Local-Linearly Convergent Algorithm}\label{sec.specific}

In this section, we prove that the Gradient-Eigenstep method, namely, \cite[Algorithm 1]{GoyeEfteBoum2024}---if the step size is sufficiently small and the penalty parameter is sufficiently large---converges linearly when initiated in the vicinity of certain second-order stationary points of the optimization problem
\bequation\label{eq.problem}
  \min_{x\in\R{n}}\ f(x)\ \ \text{subject to } (\st)\ \ c(x)=0
\eequation
when, among other properties, the functions $f$ and $c$ satisfy the following.

\bassumption\label{ass.smooth}
  The objective function $f : \R{n} \to \R{}$ and constraint function $c : \R{n} \to \R{m}$ are twice-continuously differentiable.
\eassumption

\noindent
We begin this section with some additional assumptions and definitions, then proceed with our theoretical convergence-rate analysis.

As in \cite{GoyeEfteBoum2024}, we assume the following about the constraint function.

\bassumption\label{ass.compact}
  There exists a positive real number $R \in \R{}_{>0}$ such that the set $\Ccal_R := \{x \in \R{n} : \|c(x)\| \leq R\}$ is compact.
\eassumption

Under Assumption~\ref{ass.smooth}, let us denote the objective gradient function as $\nabla f : \R{n} \to \R{n}$, the objective Hessian function as $\nabla^2 f : \R{n} \to \R{n \times n}$, the constraint derivative function as $\nabla c : \R{n} \to \R{n \times m}$, and for each $j \in [m]$ the $j$th constraint Hessian function as $\nabla^2 [c]_j : \R{n} \to \R{n \times n}$.  Here and throughout, for a function $c$, we use $[c]_j$ to denote its $j$th component.  Similarly as in~\cite{GoyeEfteBoum2024}, we make the following assumption about the constraint Jacobian function, i.e., $\nabla c^T$.  Here, the margin $\delta$ appears in our analysis as a tolerance to ensure that trial steps computed by an algorithm remain within the set $\hat\Ccal_R$.

\bassumption\label{ass.licq}
  There exists $\delta \in \R{}_{>0}$, a bounded open set $\hat\Ccal_R$ containing $\Ccal_R + \{v \in \R{n} : \|v\| \leq \delta\}$ and a positive real number $\sigma_{\min} \in \R{}_{>0}$ such that, for any $x \in \hat\Ccal_R$, the constraint Jacobian has $\sigma_m(\nabla c(x)^T) \geq \sigma_{\min}$.
\eassumption

The Lagrangian of \eqref{eq.problem} is $L:\R{n} \times \R{m} \to \R{}$ defined by
\bequation\label{eq.Lagrangian}
  L(x,\yhat) = f(x) + c(x)^T\yhat,
\eequation
where $\yhat \in \R{m}$ is known as the Lagrange multiplier.  Defined with respect to~$L$, an approximate second-order stationary point of~\eqref{eq.problem} is the following.

\begin{definition}\label{def.approximate.stationary}
  For a pair of real numbers $(\epsilon, \zeta) \in \R{}_{\geq0} \times \R{}_{\geq0}$, a point $x\in\R{n}$ is $(\epsilon, \zeta)$-stationary for \eqref{eq.problem} if there exists $\yhat \in \R{m}$ such that
  \begin{align*}
    \|\nabla L(x, \yhat)\| &\leq \epsilon \\ \text{and}\ \ 
    d^T \nabla^2_{xx} L(x,\yhat)d &\geq -\zeta \|d\|^2\ \ \text{for all}\ \ d \in \Null(\nabla c(x)^T).
  \end{align*}
\end{definition}

Due to the first of the second-order stationarity conditions in Definition~\ref{def.approximate.stationary}, for any $x \in \R{n}$ a particular Lagrange multiplier of interest is any that minimizes a norm of the gradient of the Lagrangian.  Such a Lagrange multiplier is also of interest since it is used within the algorithm that we analyze.  Let us define an $\ell_2$-norm least-squares Lagrange multiplier function $y : \hat\Ccal_R \to \R{m}$ by
\bequation\label{eq.y.lse.formulation}
  y(x) = \arg\min_{\yhat \in \R{m}} \|\nabla f(x) + \nabla c(x) \yhat\|^2,
\eequation
and note that, under Assumption~\ref{ass.licq}, for any $x \in \hat\Ccal_R$ the unique solution is
\bequation\label{eq.y}
  y(x) = -(\nabla c(x)^T \nabla c(x))^{-1} \nabla c(x)^T \nabla f(x) = -\nabla c(x)^\dag \nabla f(x).
\eequation

We now introduce Fletcher's augmented Lagrangian function.

\begin{definition}\label{def.fal}
  Using $y$ defined by \eqref{eq.y}, Fletcher's augmented Lagrangian with penalty parameter $\rho \in \R{}_{>0}$ is the function $F_\rho : \hat\Ccal_R \to \R{}$ defined by 
  \bequationNN
    F_\rho(x) = f(x) + c(x)^Ty(x) + \rho \|c(x)\|^2.
  \eequationNN
\end{definition}

We also include the following additional assumption about the Lagrangian with respect to the least-squares Lagrange multiplier function and the constraint violation measure $\|c\|^2$, i.e., the two component functions defining Fletcher's augmented Lagrangian for any penalty parameter.

\bassumption\label{ass.Lagrangian}
  Along with Assumptions~\ref{ass.smooth}, \ref{ass.compact}, and \ref{ass.licq}, the function $y$ defined by \eqref{eq.y} is twice-continuously differentiable over $\hat\Ccal_R$, and over $\hat\Ccal_R$ the functions $f + c^Ty$ and $\|c\|^2$ have Hessian functions that are Lipschitz continuous in the sense that there exists $(M_f, M_c) \in \R{}_{>0} \times \R{}_{>0}$ such that, for all $(x,\xbar) \in \hat\Ccal_R \times \hat\Ccal_R$, one has that
  \begin{align*}
    \|\nabla^2 (f + c^Ty)(x) - \nabla^2 (f + c^Ty)(\xbar)\| &\leq M_f \|x - \xbar\| \\ \text{and}\ \ 
    \|\nabla^2 (\|c\|^2)(x) - \nabla^2 (\|c\|^2)(\xbar)\| &\leq M_c \|x - \xbar\|.
  \end{align*}
  Under these conditions, for any $\rho \in \R{}_{>0}$, define $M_\rho := M_f + \rho M_c$.
\eassumption

We now commence our analysis.  For our first lemma, we introduce a uniform bound for the norms of various quantities and state useful properties of~$F_\rho$ at approximate second-order stationary points.

\begin{lemma}\label{lemma.F.psd}
  Under Assumptions \ref{ass.smooth}, \ref{ass.compact}, \ref{ass.licq}, and \ref{ass.Lagrangian} there is $\kappa \in \R{}_{>0}$ with
  \begin{multline}\label{ineq.bound.bar.kappa}
    \max_{x \in \Ccal_R} \bigg\{ \|\nabla f(x)\|, \|y(x)\|, \|\nabla y(x)\|, \|\nabla c(x)\|, \\
    \|\nabla^2 f(x)\|, \max_{j\in[m]} \left\{ \|\nabla^2 [y]_j(x)\|, \|\nabla^2 [c]_j(x)\| \right\} \bigg\}\le\kappa,
  \end{multline}
  and along with any $\rho \in \R{}_{>0}$ one has
  \bequation\label{sigma.max.hess.F}
    \lambda_{R,\rho} := \sup_{x \in \Ccal_R} \{\lambda_1(\nabla^2 F_\rho(x))\} < \infty.
  \eequation
  Moreover, for any $\epsilon \in (0,R]$ and $\beta \in \R{}_{>0}$, if $x \in \R{n}$ has
  \bequation\label{eq.epsbeta}
    \baligned
      \|\nabla L(x,y(x))\| &\leq \epsilon \\ \text{and}\ \ 
      d^T\nabla_{xx}^2 L(x,y(x)) d &\geq \beta \|d\|^2\ \ \text{for all}\ \ d \in \Null(\nabla c(x)^T),
    \ealigned
  \eequation
  then for any $\rho \in \R{}_{>0}$ and with positive real numbers
  \bequation\label{eq.def.eta_1,2}
    \baligned
      \omega_\rho &:= 1 + \kappa + 2\rho\kappa \\ \text{and}\ \ 
      \eta_\rho &:= \frac{2\sqrt{m}\kappa}{\sigma_{\min}} + m\kappa + 2 m \rho \kappa
    \ealigned
  \eequation
  one finds that
  \bequation\label{ineq.fletcher.condition}
    \baligned
      \|c(x)\| &\leq R, \\
      \|\nabla F_\rho(x)\| &\leq \omega_\rho \epsilon, \\ \text{and}\ \ 
      \lambda_n(\nabla^2 F_\rho(x)) &\geq \min \{\beta - \eta_\rho \epsilon, 2 \rho \sigma_{\min}^2 - (\kappa + m \kappa^2) - \eta_\rho \epsilon \}.
    \ealigned
  \eequation
\end{lemma}
\bproof
  Consider~\eqref{ineq.bound.bar.kappa}. Continuity of $\nabla f(x)$, $\nabla c(x)$, $\nabla^2 f(x)$, and $\nabla^2 [c]_j(x)$ for all $j\in[m]$ follow by Assumption~\ref{ass.smooth}. Continuity of $y(x)$, $\nabla y(x)$, and $\nabla^2 [y]_j(x)$ for all $j\in[m]$ follow by Assumptions~\ref{ass.smooth}, \ref{ass.licq}, and \ref{ass.Lagrangian}. Since $\Ccal_R$ is compact by Assumption~\ref{ass.compact}, the existence of $\kappa \in \R{}_{>0}$ follows from continuity of these functions, continuity of norms, the fact that compositions of continuous functions are continuous, and the Weierstrauss Extreme Value Theorem.

  Consider next~\eqref{sigma.max.hess.F}.  The fact that $F_\rho$ is twice-continuously differentiable follows under Assumption~\ref{ass.Lagrangian}, and by continuity of eigenvalues of $\nabla^2 F_\rho$ and compactness of $\Ccal_R$ it follows that the supremum in \eqref{sigma.max.hess.F} is finite.

  Consider now~\eqref{ineq.fletcher.condition}. Consider arbitrary $\epsilon \in (0,R]$, $\beta \in \R{}_{>0}$, $x \in \R{n}$ satisfying \eqref{eq.epsbeta}, and $\rho \in \R{}_{>0}$. Therefore,
  \bequation\label{ineq.nablax.nablay}
    \max\{\|\nabla f(x) + \nabla c (x) y (x)\|, \|c (x)\|\} \leq  \|\nabla L (x, y(x))\| \leq \epsilon,
  \eequation
  which with $\epsilon \leq R$ gives the first inequality in \eqref{ineq.fletcher.condition}.  Thus, $x \in \Ccal_R$, so the triangle inequality and submultiplicity of the matrix 2-norm give
  \begin{align}
    \|\nabla F_\rho(x)\|
      &\leq \|\nabla f(x) + \nabla c(x) y(x)\| + (\|\nabla y(x)\| + 2 \rho \|\nabla c(x)\|) \|c(x)\| \nonumber \\
      &\leq \|\nabla f(x) + \nabla c(x) y(x)\| + (\kappa + 2 \rho \kappa) \|c(x)\|. \label{eq.fletcher.al.grad}
  \end{align}
  Combining \eqref{eq.fletcher.al.grad} with \eqref{ineq.nablax.nablay}, one obtains the second inequality in \eqref{ineq.fletcher.condition}.  Our final aim is to prove the third inequality in \eqref{ineq.fletcher.condition}.  Toward this end, note that
  \begin{align}
    &\ \nabla^2 F_\rho(x) \nonumber \\
    =&\ \nabla^2_{xx} L(x, y(x)) + \nabla y(x) \nabla c(x)^T + \nabla c(x) \nabla y(x)^T \nonumber \\
    &\ + 2\rho \nabla c(x) \nabla c(x)^T + \sum_{j\in[m]} (\nabla^2 [y]_j(x) + 2 \rho \nabla^2 [c]_j(x)) [c]_j(x). \label{eq.hess.F}
  \end{align}
  Let us now express this Hessian in an equivalent form involving a decomposition into orthogonal spaces defined by the constraint Jacobian at $x$.  Let us first derive an expression for $\nabla y(x)\nabla c(x)^T$ by differentiating the linear system that defines $y(x)$ through \eqref{eq.y}.  Specifically, by \eqref{eq.y}, one finds that
  \bequation\label{eq.grad.y}
    \nabla c (x)^T\nabla c (x) y (x) = -\nabla c (x)^T \nabla f(x).
  \eequation
  Abbreviating notation and differentiating the left-hand side yields
  \begin{align*}
    &\ \nabla(\nabla c^T\nabla c y) |_x \\
    =&\ \nabla ( \nabla c^T|_x \nabla c|_x y) |_x + \nabla (\nabla c^T \nabla c|_x y|_x ) |_x + \nabla (\nabla c^T|_x \nabla c\ y|_x) |_x \\
    =&\ \nabla y|_x \nabla c^T |_x \nabla c|_x + \bbmatrix \nabla^2 [c]_1|_x \nabla c|_x y|_x & \cdots
    & \nabla^2 [c]_m|_x \nabla c|_x y|_x
    \ebmatrix \\
    &\ + \(\sum_{j\in[m]} \nabla^2 [c]_j|_x [y]_j|_x \) \nabla c|_x,
  \end{align*}
  while at the same time differentiating the right-hand side yields
  \begin{align*}
    \nabla ( -\nabla c^T\nabla f) |_x
    =&\ -\nabla (\nabla c^T|_x \nabla f)|_x - \nabla (\nabla c^T\nabla f|_x) |_x \\
    =&\ -\nabla^2 f|_x \nabla c|_x -\bbmatrix \nabla^2 [c]_1|_x \nabla f|_x & \cdots & \nabla^2 [c]_m|_x \nabla f|_x \ebmatrix.
  \end{align*}
  Combining these derivations and rearranging yields
  \begin{align}
    &\ \nabla y(x) \nabla c(x)^T \nabla c(x) \nonumber \\
    =&\ -\nabla^2_{xx} L(x,y(x)) \nabla c(x) \nonumber \\
    &\ - \underbrace{\bbmatrix \nabla^2 [c]_1(x) \nabla_x L(x,y(x)) & \cdots & \nabla^2 [c]_m(x) \nabla_x L(x,y(x)) \ebmatrix}_{=: \Ecal(x)}. \label{ineq.important}
  \end{align}
  Multiplying \eqref{ineq.important} on the right by $(\nabla c(x)^T \nabla c(x))^{-1} \nabla c(x)^T$ yields
  \begin{align}
    &\ \nabla y(x) \nabla c(x)^T \nonumber \\
    =&\ -(\nabla^2_{xx} L(x,y(x)) \nabla c(x) + \Ecal(x)) (\nabla c(x)^T \nabla c(x))^{-1} \nabla c(x)^T \nonumber \\
    =&\ -\nabla^2_{xx} L(x,y(x)) \Rcal(\nabla c(x)) -\Ecal(x) \nabla c(x)^\dag. \label{eq.hess.e}
  \end{align}
  On the other hand, denoting $\Rcal(x) := \Rcal(\nabla c(x))$, $\Ncal(x) := \Ncal(\nabla c(x)^T)$, and $H_{xx} := \nabla^2_{xx} L(x,y(x))$ and using $\Rcal(x) + \Ncal(x) = I$, one finds that
  \begin{align}
    H_{xx} - H_{xx}\Rcal(x) - \Rcal (x)H_{xx}
    &= H_{xx} \Ncal(x) - \Rcal(x) H_{xx} \nonumber \\
    &=(\Ncal(x) + \Rcal(x)) H_{xx} \Ncal(x) - \Rcal(x) H_{xx} \nonumber \\
    &= \Ncal(x) H_{xx} \Ncal(x) - \Rcal(x) H_{xx} \Rcal(x). \label{eq.hess.simp}
  \end{align}
  Combining \eqref{eq.hess.F}, \eqref{eq.hess.e}, and \eqref{eq.hess.simp} now yields the expression
  \begin{align}
    \nabla^2 F_\rho(x)
    =&\ \Ncal(x) H_{xx} \Ncal(x) -\Rcal(x) H_{xx} \Rcal(x) + 2 \rho \nabla c(x) \nabla c(x)^T \nonumber \\
    &\ -\Ecal(x) \nabla c(x)^\dag -(\nabla c(x)^\dag)^T \Ecal(x)^T \nonumber \\
    &\ + \sum_{j\in[m]} (\nabla^2 [y]_j(x) + 2 \rho \nabla^2 [c]_j(x)) [c]_j(x). \label{eq.fl.hess.rewrite}
  \end{align}
  Let us now observe that from the definition of $\Ecal(x)$, norm inequalities, submultiplicity of the matrix 2-norm, $x \in \Ccal_R$, \eqref{ineq.bound.bar.kappa}, and \eqref{ineq.nablax.nablay} that
  \begin{align}
    \|\Ecal(x)\|
    &= \|\Ecal(x)^T\| \nonumber \\
    &= \max_{d\in\R{n} \st \|d\|=1} \|\Ecal(x)^Td\| \nonumber \\
    &= \max_{d\in\R{n} \st \|d\|=1} \left\| \bbmatrix \nabla_x L(x,y(x))^T \nabla^2 [c]_1(x)^Td \\ \vdots \\ \nabla_x L(x,y(x))^T \nabla^2 [c]_m(x)^Td \ebmatrix \right\| \nonumber \\
    &= \max_{d\in\R{n} \st \|d\|=1} \sqrt{\sum_{j\in[m]} \( \nabla_x L(x,y(x))^T \nabla^2 [c]_j(x)^Td\)^2} \nonumber \\
    &\leq \sqrt{m} \max_{d\in\R{n} \st \|d\|=1} \(\max_{j\in[m]} |\nabla_x L(x,y(x))^T \nabla^2 [c]_j(x)^Td| \) \nonumber \\
    &\leq \sqrt{m} \max_{d\in\R{n} \st \|d\|=1} \(\max_{j\in[m]} \|\nabla_x L(x,y(x))\| \|\nabla^2 [c]_j(x)\| \|d\|\) \nonumber \\
    &\leq \sqrt{m} \max_{d\in\R{n} \st \|d\|=1} (\epsilon\kappa\|d\|) = \sqrt{m} \kappa \epsilon. \label{ineq.norm.Ecal}
  \end{align}
  Now consider arbitrary $d \in \R{n} \setminus \{0\}$ decomposed as $d \equiv d_\Ncal + d_\Rcal$, where $d_\Ncal \in \Null(\nabla c(x)^T)$ and $d_\Rcal \in \Range(\nabla c(x))$.  One has from \eqref{eq.fl.hess.rewrite} that
  \begin{align}
    &\ d^T\nabla^2 F_\rho(x) d \nonumber \\
    =&\ d^T\Ncal(x) \nabla^2_{xx} L(x,y(x)) \Ncal(x) d - d^T \Rcal(x) \nabla^2_{xx} L(x,y(x)) \Rcal(x) d \nonumber \\
    &\ + 2\rho d^T \nabla c(x)\nabla c(x)^T d - 2 d^T \Ecal(x) \nabla c(x)^\dag d \nonumber \\
    &\ + d^T\(\sum_{j\in[m]} \(\nabla^2 [y]_j(x) + 2\rho \nabla^2 [c]_j(x)\) [c]_j(x)\)d \nonumber \\
    =&\ d_\Ncal^T \nabla^2_{xx} L(x,y(x)) d_\Ncal - d_\Rcal^T \nabla^2_{xx} L(x,y(x)) d_\Rcal \nonumber \\
    &\ + 2\rho \|\nabla c(x)^T d_\Rcal\|^2 - 2 d^T\Ecal(x) \nabla c(x)^\dag d \nonumber \\
    &\ + d^T \(\sum_{j\in[m]} \(\nabla^2 [y]_j(x) + 2\rho \nabla^2 [c]_j(x)\) [c]_j(x)\) d, \label{ineq.hess.dd}
  \end{align}
  where for the latter two terms one has from \eqref{ineq.bound.bar.kappa}, \eqref{ineq.nablax.nablay}, \eqref{ineq.norm.Ecal}, and the fact that the spectral norm of the pseudoinverse of a matrix is at most the reciprocal of the matrix's smallest singular value \cite[Chapter 21]{gallier2019linear} that
  \begin{align}
    &\ -2d^T \Ecal(x) \nabla c(x)^\dag d + d^T \(\sum_{j\in[m]} (\nabla^2 [y]_j(x) + 2\rho\nabla^2 [c]_j(x)) [c]_j(x)\)d \nonumber \\
    \geq&\ -2 \|\Ecal(x)\| \|\nabla c(x)^\dag\| \|d\|^2 - \left\| \sum_{j\in[m]} (\nabla^2 [y]_j(x) + 2\rho\nabla^2 [c]_j(x)) [ c]_j(x) \right\| \|d\|^2 \nonumber \\
    \geq&\ -\frac{2\sqrt{m}\kappa}{\sigma_{\min}} \epsilon \|d\|^2 -\sum_{j\in[m]} ( \|\nabla^2 [y]_j(x)\| + 2\rho\|\nabla^2 [c]_j(x)\|) |[c]_j(x)| \|d\|^2 \nonumber \\ 
    \geq&\ -\frac{2\sqrt{m}\kappa}{\sigma_{\min}} \epsilon \|d\|^2 -\sum_{j\in[m]} ( \|\nabla^2 [y]_j(x)\| + 2\rho \|\nabla^2 [c]_j(x)\| ) \|c(x)\| \|d\|^2 \nonumber \\
    \geq&\ -\frac{2\sqrt{m}\kappa}{\sigma_{\min}} \epsilon \|d\|^2 - (m \kappa + 2 m \rho\kappa) \epsilon \|d\|^2. \label{ineq.bound.sum}
  \end{align}
  Since $d_\Ncal \in \Null(\nabla c(x)^T)$, Definition~\ref{def.approximate.stationary} gives $d_\Ncal^T \nabla_{xx}^2 L(x,y(x)) d_\Ncal \geq \beta \|d_\Ncal\|^2$.  With \eqref{ineq.bound.bar.kappa}, the triangle inequality, and submultiplicity of the matrix 2-norm,
  \begin{align}
    \|\nabla^2_{xx} L(x,y(x))\|
      &= \left\|\nabla^2f(x) + \sum_{j\in[m]} \nabla^2 [c]_j(x)[y]_j(x) \right\| \nonumber \\
      &\leq \|\nabla^2 f(x)\| + \sum_{j\in[m]} \|\nabla^2 [c]_j(x)\| \|y(x)\| \leq \kappa + m \kappa^2. \label{ineq.hess.lag}
  \end{align}
  At the same time, let us derive a lower bound for $\|\nabla c(x)^Td_\Rcal\|$.  Let $\nabla c(x)$ have the singular value decomposition $\sum_{i\in[m]} u_i \sigma_i v_i^T$ where $(u_i,\sigma_i,v_i) \in \R{n} \times \R{} \times \R{m}$ for all $i\in[m]$ with $\{u_i\}_{i\in[m]}$ and $\{v_i\}_{i\in[m]}$ being sets of orthonormal vectors. Then, by Assumption~\ref{ass.licq} and the fact that $d_\Rcal\in\Range(\nabla c(x))$,
  \begin{align}
    \|\nabla c(x)^Td_\Rcal\|^2
    &= \left\| \sum_{i\in[m]} v_i \sigma_i u_i^T d_\Rcal \right\|^2 = \sum_{i\in[m]} \|v_i \sigma_i u_i^T d_\Rcal\|^2 = \sum_{i\in[m]} \sigma_i^2 (u_i^T d_\Rcal)^2 \nonumber \\
&\ge\sigma_{\min}^2\sum_{i=1}^m(u_i^Td_\Rcal)^2=\sigma_{\min}^2\|d_\Rcal\|^2. \label{ineq.nabla cd}
  \end{align}
  Combining \eqref{ineq.hess.dd}, \eqref{ineq.bound.sum}, \eqref{ineq.hess.lag}, \eqref{ineq.nabla cd}, and the decomposition $d = d_\Ncal + d_\Rcal$,
  \begin{align}
    &\ d^T\nabla^2 F_\rho(x) d \nonumber \\
    \geq&\ \beta \|d_\Ncal\|^2 + (2\rho\sigma_{\min}^2 - (\kappa+ m \kappa^2 )) \|d_\Rcal\|^2 \nonumber \\
    &\ - \(\frac{2 \sqrt{m} \kappa}{\sigma_{\min}} + m \kappa + 2m \rho \kappa\) \epsilon \|d\|^2 \nonumber \\
    =&\ \(\beta - \(\frac{2\sqrt{m} \kappa}{\sigma_{\min}} + m \kappa + 2 m \rho \kappa \) \epsilon\) \|d_\Ncal\|^2 \nonumber \\
    &\ + \(2\rho\sigma_{\min}^2 - (\kappa + m \kappa^2) - \(\frac{2\sqrt{m}\kappa}{\sigma_{\min}} + m \kappa + 2 m \rho \kappa \) \epsilon\) \|d_\Rcal\|^2. \label{ineq.F.Hessian.cur}
  \end{align}
  Combined with the definition of $\eta_\rho$ in \eqref{eq.def.eta_1,2}, the proof is complete.
  \qed
\eproof

We make a few observations regarding Lemma~\ref{lemma.F.psd}.  First, the bounds on the norms of $c(x)$ and $\nabla F_\rho(x)$ in \eqref{ineq.fletcher.condition} follow from the fact that~$x$ is assumed to satisfy \eqref{eq.epsbeta} with $\epsilon \leq R$; the factor on the right-hand side of the second inequality in \eqref{ineq.fletcher.condition} merely accounts for the difference between the gradients of the Lagrangian and Fletcher's augmented Lagrangian.  Second, we note that our lower bound on the smallest eigenvalue of the Hessian of Fletcher's augmented Lagrangian in~\eqref{ineq.fletcher.condition} can be related to \cite[Theorem 6]{estrin2020}, which gives a lower bound for $\rho$ to ensure positive semidefiniteness of $\nabla^2 F_\rho(x)$, namely,
\bequation\label{eq.estrin}
  \rho \geq \max\{0, \lambda_1(\nabla c(x)^\dag\nabla^2_{xx} L(x) (\nabla c(x)^\dag)^T)\}.
\eequation
In particular, \eqref{ineq.F.Hessian.cur} with $\epsilon=0$ yields
\bequationNN
  d^T \nabla^2 F_\rho(x) d \geq \beta \|d_\Ncal\|^2 + (2\rho\sigma_{\min}^2 - (\kappa + m \kappa^2)) \|d_\Rcal\|^2,
\eequationNN
so positive semidefiniteness of $\nabla^2 F_\rho(x)$ follows from $\rho \geq \frac{\kappa + m \kappa^2}{2\sigma_{\min}^2}$.  This is of the same form as \eqref{eq.estrin} since, as mentioned in the proof of the previous lemma, $\|\nabla c(x)^\dag\| \leq \frac{1}{\sigma_{\min}}$ and from the definition of $\kappa$ one has $\|\nabla^2_{xx} L(x)\| \leq \kappa + m\kappa^2$.  Finally, let us make two further observations related to \eqref{ineq.F.Hessian.cur}.  That is, for directions in $\Null(\nabla c(x)^T)$, a small $\epsilon$ relative to $\beta$ is sufficient to ensure positive curvature with respect to $\nabla^2 F_\rho(x)$. However, for directions in $\Range(\nabla c(x))$, the combination of a large $\rho$ and a small $\epsilon$ are needed to ensure positive curvature with respect to $\nabla^2 F_\rho(x)$.  This occurs since \eqref{eq.epsbeta} corresponds to sufficiently positive curvature of $\nabla_{xx}^2 L(x,y(x))$ in $\Null(\nabla c(x)^T)$, but to ensure positive curvature of $\nabla^2 F_\rho(x)$ in $\Range(\nabla c(x))$ one needs a large penalty parameter~$\rho$ to exploit the curvature of the Hessian of $\|c(x)\|^2$.

We now prove our main theorem of this section.

\begin{theorem}\label{theo.fl.strongly.convex}
  Suppose that Assumptions \ref{ass.smooth}, \ref{ass.compact}, \ref{ass.licq}, and \ref{ass.Lagrangian} hold and let $\kappa \in \R{}_{>0}$ satisfy \eqref{ineq.bound.bar.kappa}. In addition, for some $\beta \in \R{}_{>0}$, suppose that
  \bequation\label{eq.rho_lower}
    \rho > \rho_\beta := \max\left\{\frac{\max\{1, \kappa\}}{\sigma_{\min}}, \frac{\beta + \kappa + m \kappa^2}{2 \sigma_{\min}^2}\right\},
  \eequation
  and that $x \in \R{n}$ satisfies \eqref{eq.epsbeta}, where with $M_\rho \in \R{}_{>0}$ defined in Assumption~\ref{ass.Lagrangian}, $(\omega_\rho,\eta_\rho) \in \R{}_{>0} \times \R{}_{>0}$ defined in \eqref{eq.def.eta_1,2}, and $\lambda_\rho \in \R{}_{>0}$ satisfying
  \bequationNN
    \lambda_\rho \geq \beta + \lambda_{R,\rho}\ \ \text{$($where $\lambda_{R,\rho}$ is defined in \eqref{sigma.max.hess.F}$)$}
  \eequationNN
  the tolerance $\epsilon \in \R{}_{>0}$ satisfies $($with $\delta$ defined in Assumption~\ref{ass.licq}$)$
  \bequation\label{ineq.requirement.theorem}
    \epsilon \leq \epsilon_\rho := \min\left\{\frac{\beta}{2\(\eta_\rho + \frac{4 M_\rho \omega_\rho}{\beta} +\(1+\frac{2}{\sigma_{\min}}\) \frac{\kappa}{\sqrt{5}}\)}, \frac{R}{\omega_\rho}, \frac{\beta\lambda_\rho}{2M_\rho\omega_{\rho}}, \frac{\lambda_\rho \delta}{2 \omega_\rho} \right\}.
  \eequation
  Furthermore, suppose that $t \in \R{}_{>0}$ satisfies
  \bequation\label{eq.stepsize.statement}
    0 < t \leq t_\rho := \min \left\{\frac{1}{\lambda_\rho}, \frac{1}{\beta - \eta_\rho\epsilon - \frac{4 M_\rho \omega_\rho}{\beta} \epsilon} \right\}.
  \eequation
  Then, for any $\epsilon' \in (0,\epsilon)$ and with the parameters
  \bequation\label{goyens.alg.parameter}
    \underbrace{(\epsilon_1, \alpha_{01}, c_1)}_{\text{(notation in \cite{GoyeEfteBoum2024})}} \equiv \(\frac{\epsilon'}{\sqrt{5}}, t, \half\),
  \eequation
  the method \cite[Algorithm 1]{GoyeEfteBoum2024} reduces to applying gradient descent with the constant step size $t$ to minimize $F_\rho$ with initial point~$x$.  Thus, in at most
  \bequation\label{eq.T}
    T = \left \lceil \log_{\frac{4}{4 - \beta t}} \frac{\sqrt{5} \omega_\rho \epsilon}{\epsilon'} \right \rceil\ \text{iterations $($where $4 - \beta t > 0$$)$}
  \eequation
  the method gives a point satisfying \eqref{eq.epsbeta} with $(\epsilon,\beta)$ replaced by $(\epsilon',\beta')$, where
  \bequation\label{eq.beta'}
      \beta' := \beta - \eta_\rho \epsilon - \frac{4 M_\rho \omega_\rho}{\beta} \epsilon - \(1 + \frac{2}{\sigma_{\min}} \) \frac{\kappa}{\sqrt{5}} \epsilon' \geq \half \beta.
  \eequation
\end{theorem}
\bproof
  Let us begin by showing that \eqref{ineq.requirement.theorem} implies \eqref{eq.beta'} as well as positivity of other quantities that are used throughout the proof, including the upper bound for the step size in \eqref{eq.stepsize.statement}.  By \eqref{ineq.requirement.theorem} and $\epsilon' \in (0,\epsilon)$, one finds
  \begin{align*}
    \beta'
      &= \beta - \eta_\rho \epsilon - \frac{4 M_\rho \omega_\rho}{\beta} \epsilon - \(1 + \frac{2}{\sigma_{\min}} \) \frac{\kappa}{\sqrt{5}} \epsilon' \\
      &\geq \beta - \(\eta_\rho + \frac{4 M_\rho \omega_\rho}{\beta} + \(1 + \frac{2}{\sigma_{\min}} \) \frac{\kappa}{\sqrt{5}}\) \epsilon \geq \half \beta,
  \end{align*}
  which shows that \eqref{eq.beta'} holds.  From these inequalities, one also finds that
  \bequation\label{eq.beta.eps.rho}
    \beta_{\rho,\epsilon} := \beta - \eta_\rho\epsilon \geq \half \beta\ \ \text{and}\ \ \tilde{\beta}_{\rho,\epsilon} := \beta - \eta_\rho\epsilon - \frac{4 M_\rho \omega_\rho}{\beta} \epsilon \geq \half \beta,
  \eequation
  which shows that the upper bound in \eqref{eq.stepsize.statement} is positive.
  
  Let us now consider gradient descent applied to minimize $F_\rho$ from $x$ with step size~$t$, where $\rho \in \R{}_{>0}$ satisfies \eqref{eq.rho_lower} and $t \in \R{}_{>0}$ satisfies \eqref{eq.stepsize.statement}, which by combining \eqref{ineq.requirement.theorem} and \eqref{eq.beta.eps.rho} can be seen to satisfy
  \bequation\label{eq.stepsize}
    0 < t \leq \frac{1}{\lambda_\rho } \leq \frac{\beta}{2 M_\rho \omega_\rho \epsilon} \leq \frac{\tilde\beta_{\rho,\epsilon}}{M_\rho \omega_\rho \epsilon}\ \text{and}\ 0 < t \leq \frac{1}{\tilde\beta_{\rho,\epsilon}}.
  \eequation
  Once this is completed, we confirm that \cite[Algorithm 1]{GoyeEfteBoum2024} with the parameters in \eqref{goyens.alg.parameter} indeed reduces to such an application of gradient descent.  Let $x^0 \gets x$ and consider the iterative method defined for all $i = 0,1,2,\dots$ by
  \bequation\label{eq.fl.gd}
    x^{i+1} \gets\ x^i - t \nabla F_\rho(x^i)\ \ \text{until}\ \ \|\nabla F_\rho(x^{i+1})\| \leq \frac{\epsilon'}{\sqrt{5}}.
  \eequation
  Our first aim is to show that $x^0 \in \Ccal_R$, $\|\nabla F_\rho(x^0)\| \leq \omega_\rho \epsilon$, and $\lambda_n(\nabla^2 F_\rho(x^0)) \geq \beta_{\rho,\epsilon}$, and that, for any $i \in \N{}$ with $i \geq 1$, one finds that
  \bequation\label{ineq.induct.end}
    \baligned
      x^i &\in \Ccal_R, \\
      \|\nabla F_\rho(x^i)\| &\leq \(1 - \half  \tilde{\beta}_{\rho,\epsilon} t\)^i \|\nabla F_\rho(x^{0})\| \\ \text{and}\ \ 
      \lambda_n(\nabla^2 F_\rho(x^{i})) &\geq \beta_{\rho,\epsilon} - M_\rho t \|\nabla F_\rho(x^{0})\| \sum_{l=0}^{i-1} \(1 - \half \tilde{\beta}_{\rho,\epsilon} t \)^l \geq \tilde{\beta}_{\rho,\epsilon},
    \ealigned
  \eequation
  where $1 - \half  \tilde{\beta}_{\rho,\epsilon} t > 0$ by \eqref{eq.stepsize}.  First, for $i=0$, by \eqref{eq.rho_lower} one has that $\rho \geq \frac{\beta + \kappa + m \kappa^2}{2\sigma_{\min}^2}$, so with \eqref{eq.beta.eps.rho} one has
  \bequation\label{ineq.rho.1}
    2 \sigma_{\min}^2 \rho - (\kappa + m \kappa^2) - \eta_\rho \epsilon \geq \beta - \eta_\rho \epsilon = \beta_{\rho,\epsilon}. 
  \eequation
  Since $x^0 = x$ satisfies \eqref{eq.epsbeta}, it follows from \eqref{ineq.requirement.theorem} and $\omega_\rho \geq 1$ that $\|c(x^0)\| \leq R$, which means $x^0 \in \Ccal_R$, and it follows from Lemma~\ref{lemma.F.psd}, \eqref{ineq.requirement.theorem}, and \eqref{ineq.rho.1} that
  \bequation\label{ineq.i-0}
    \begin{aligned}
      \|\nabla F_\rho(x^0)\| &\leq \omega_\rho \epsilon \\ \text{and}\ \  
      \lambda_n(\nabla^2 F_\rho(x^0)) &\geq \min\{\beta - \eta_\rho \epsilon, 2\sigma_{\min}^2 \rho - (\kappa + m \kappa^2) - \eta_\rho \epsilon\} \\
      &\geq \beta_{\rho,\epsilon},
    \end{aligned}
  \eequation
  as claimed.  Now consider arbitrary $i \in \N{}$ with $i \geq 1$ such that \eqref{ineq.induct.end} holds. To complete the induction, our aim is to prove that all of~\eqref{ineq.induct.end} holds with $i$ replaced by $i+1$ as well.  Denote $G_i := \nabla F_\rho(x^i)$ for all $i \geq 0$.  Furthermore, note that, by \eqref{ineq.requirement.theorem} and the fact that $t \leq 1/\lambda_\rho$, one has
  \bequationNN
    \|x^{i+1} - x^i\| = t\|G_i\| \leq t\|G_0\| \leq t \omega_\rho \epsilon \leq \frac{\omega_\rho \epsilon}{\lambda_\rho} \leq \frac{\delta}{2},
  \eequationNN
  meaning that the line segment $[x^i,x^{i+1}]$ lies in $\hat\Ccal_R$ and Assumption~\ref{ass.Lagrangian} applies along the segment.  Next, one has from Assumptions~\ref{ass.smooth} and \ref{ass.Lagrangian} that
  \begin{align}
    \|G_{i+1}\|
      &\leq \left\|G_i + \nabla^2 F_\rho(x^i) (x^{i+1} - x^i) \right\| + \half  M_\rho \|x^{i+1}-x^{i}\|^2 \nonumber \\
      &= \|G_i - t \nabla^2 F_\rho(x^i)G_i\| + \half  M_\rho t^2 \|G_i\|^2 \nonumber \\
      &\leq \|I - t \nabla^2 F_\rho(x^i)\| \|G_i\| + \half  M_\rho t^2 \|G_i\|^2. \label{ineq.taylor.exp}
  \end{align}
  By $x^i \in \Ccal_R$ from \eqref{ineq.induct.end}, \eqref{sigma.max.hess.F}, $\lambda_\rho \geq \lambda_{R,\rho}$, and \eqref{eq.stepsize}, one has $t \lambda_1(\nabla^2 F_\rho(x^i)) \leq t \lambda_\rho  \leq 1$, which with \eqref{ineq.induct.end} gives
  \begin{align}
    &\ \|I - t \nabla^2 F_\rho(x^i)\| \nonumber \\
    =&\ \max_{v\in\R{n}\st\|v\|=1} \left|v^Tv-t v^T\nabla^2  F_\rho(x^i)v\right| \nonumber \\
    =&\ \max\left\{ |1 - t \lambda_n(\nabla^2 F_\rho(x^i))|, |1 - t \lambda_1(\nabla^2 F_\rho(x^i)) | \right\} \nonumber \\
    =&\ 1 - t \lambda_n( \nabla^2  F_\rho(x^i)) \nonumber \\
    \leq&\ 1 - t \tilde{\beta}_{\rho,\epsilon}. \label{ineq.bound.1}
  \end{align}
  At the same time, combining \eqref{eq.stepsize}, \eqref{ineq.induct.end}, and \eqref{ineq.i-0} yields
  \bequation\label{ineq.g0}
    M_\rho t \|G_i\| \leq  M_\rho t \|G_0\| \leq M_\rho t \omega_\rho \epsilon \le \tilde{\beta}_{\rho,\epsilon}.
  \eequation
  Now combining \eqref{eq.stepsize}, \eqref{ineq.taylor.exp}, \eqref{ineq.bound.1}, and \eqref{ineq.g0} one obtains that
  \begin{align}
    \|G_{i+1}\|
      &\leq \(1- t \tilde{\beta}_{\rho,\epsilon} + \half M_\rho t^2 \|G_i\|\) \|G_i\|\nonumber \\
      &\le\(1 - \half t \tilde{\beta}_{\rho,\epsilon}\)\|G_i\|. \label{ineq.G1.G0}
  \end{align}
  This gives the first inequality in~\eqref{ineq.induct.end}, and in fact shows that $\|G_{i+1}\| \leq \|G_i\|$ for all $i \in \N{}$.  Also, from \cite[Theorem 6.6]{Stew1973} (or see \cite[Eq.~(3)]{Stew1979}), one finds
  \begin{align*}
    &\ |\lambda_n(\nabla^2 F_\rho(x^{i+1})) - \lambda_n(\nabla^2 F_\rho(x^i))| \\
    \leq&\ \|\nabla^2 F_\rho(x^{i+1}) - \nabla^2 F_\rho(x^i)\| \leq M_\rho\|x^{i+1}-x^i\|= M_\rho t \|G_i\|.
  \end{align*}
  Combining \eqref{ineq.requirement.theorem}, \eqref{eq.beta.eps.rho}, \eqref{eq.stepsize}, \eqref{ineq.induct.end}, and \eqref{ineq.G1.G0}, one now obtains that
  \begin{align}
    &\ \lambda_n(\nabla^2 F_\rho(x^{i+1})) \nonumber \\
    \geq&\ \lambda_n(\nabla^2 F_\rho(x^i)) - |\lambda_n(\nabla^2 F_\rho(x^{i+1})) - \lambda_n(\nabla^2 F_\rho(x^i))| \nonumber \\
    \geq&\ \beta_{\rho,\epsilon} - M_\rho t \|G_0\| \sum_{l=0}^{i-1}\(1-\half\tilde{\beta}_{\rho,\epsilon}t \)^{l} - M_\rho t\|G_i\| \nonumber \\
    \geq&\ \beta_{\rho,\epsilon}-M_\rho t \|G_0\| \sum_{l=0}^i\(1-\half\tilde{\beta}_{\rho,\epsilon}t\)^{l} \nonumber \\
    \geq&\ \beta_{\rho,\epsilon}- M_\rho t \|G_0\| \sum_{l=0}^\infty\(1-\half\tilde{\beta}_{\rho,\epsilon}t \)^{l} \nonumber \\
    =&\ \beta_{\rho,\epsilon} - M_\rho t \|G_0\| \frac{2}{\tilde{\beta}_{\rho,\epsilon} t} \nonumber \\
    \geq &\ \beta_{\rho,\epsilon} - \frac{4 M_\rho \omega_\rho \epsilon}{\beta} = \tilde{\beta}_{\rho,\epsilon} .\label{ineq.strong.convex.2}
  \end{align}
  All that remains to complete the induction is to show that $x^{i+1}\in\Ccal_R$.  Toward this end, let us employ \cite[Corollary 2.7]{GoyeEfteBoum2024}.  Employing this result requires that the singular values of $\nabla c$ are bounded away from zero and that $\Ccal_R$ is compact, which hold here by Assumptions~\ref{ass.compact} and \ref{ass.licq}.  It also requires that
  \bequation\label{ineq.ensreu.tau}
    \rho  > \max_{x\in\Ccal_{ R}} \frac{ \max\{1,\|\nabla y (x)\|\}}{\sigma_{\min}(\nabla c(x)) }.
  \eequation
  This also holds here since for any $x \in \Ccal_R$ one has $\sigma_{\min}(\nabla c(x)) \geq \sigma_{\min}$ by Assumption~\ref{ass.licq} and $\|\nabla y(x)\| \leq \kappa$ by \eqref{ineq.bound.bar.kappa}, which when combined with \eqref{eq.rho_lower} shows that \eqref{ineq.ensreu.tau} holds. Now since \eqref{ineq.G1.G0} and \eqref{ineq.strong.convex.2} give
  \bequation\label{G.geometric}
    \|G_{i+1}\| \leq \tilde{\epsilon}_{i+1} := \(1 - \half  \tilde{\beta}_{\rho,\epsilon} t\)^{i+1} \|G_0\|\ \ \text{and}\ \ \nabla^2 F_\rho(x^{i+1}) \succeq \tilde{\beta}_{\rho,\epsilon} I,
  \eequation
  by \cite[Corollary 2.7]{GoyeEfteBoum2024} one has $x^{i+1}$ is an $(\tilde{\epsilon}_{i+1}, 2\tilde{\epsilon}_{i+1}, -\tilde{\beta}_{\rho,\epsilon}+(1+\frac{2}{\sigma_{\min}})\kappa\tilde{\epsilon}_{i+1})$ stationary point of \eqref{eq.problem} in the sense (defined in \cite{GoyeEfteBoum2024}) that 
  \begin{align}
    &\ \|c(x^{i+1})\| \leq  \tilde{\epsilon}_{i+1},\ \ \|\nabla_x L(x^{i+1},y(x^{i+1}))\| \leq 2\tilde{\epsilon}_{i+1}, \nonumber \\
    \text{and} &\ d^T \nabla_{xx}^2 L(x^{i+1},y(x^{i+1}))d \geq \( \tilde{\beta}_{\rho,\epsilon} - \(1+\frac{2}{\sigma_{\min}}\) \kappa \tilde{\epsilon}_{i+1}\) \|d\|^2 \nonumber \\
    \text{for all} &\ d \in \Null(\nabla c(x^{i+1})^T). \label{ineq.rela1}
  \end{align}
  Noting that $1 - \half \tilde{\beta}_{\rho,\epsilon} t < 1$, one has $\tilde{\epsilon}_{i+1} \leq \|G_0\|$. Combined with the fact that $\|G_0\| \leq \omega_\rho \epsilon\le R$ by \eqref{ineq.requirement.theorem}, one finds that $\|c(x^{i+1})\|\le \tilde{\epsilon}_{i+1}\le R$. Thus, $x^{i+1}\in\Ccal_R$ has been proved, and the desired induction to prove \eqref{ineq.induct.end} is complete.
  
  Let us now give an upper bound for the number of iterations until the gradient descent method defined by \eqref{eq.fl.gd} terminates.  Let $T' \in \N{}$ be the smallest positive integer such that $\tilde\epsilon_{T'} \leq \frac{\epsilon'}{\sqrt{5}}$.  Then,
  \bequation\label{eq.Tform}
    T' \leq \left\lceil \log_{1-\half \tilde{\beta}_{\rho,\epsilon} t}\frac{\epsilon'}{\sqrt{5}\|G_0\|} \right\rceil,
  \eequation
  and since $\|G_{T'}\| \leq \tilde\epsilon_{T'}$ one has $T \leq T'$.  Since the definition of $\lambda_\rho$, \eqref{eq.beta.eps.rho}, \eqref{eq.stepsize}, and \eqref{ineq.induct.end} imply that
  \bequationNN
    4 - \beta t \geq 4 - \frac{\beta}{\lambda_\rho} \geq 4 - \frac{\beta}{\lambda_n(\nabla^2 F_\rho(x^T))} \geq 4 - \frac{\beta}{\tilde\beta_{\rho,\epsilon}} \geq 2,
  \eequationNN
  one finds with \eqref{eq.beta.eps.rho} that
  \bequationNN
    \frac{4}{4-\beta t} = \frac{1}{1 - \frac14 \beta t} \leq \frac{1}{1 - \frac12 \tilde{\beta}_{\rho,\epsilon} t},
  \eequationNN
  so by considering \eqref{eq.Tform} and $\|G_0\| \leq \omega_\rho \epsilon$ one finds that
  \begin{align*}
    \log_{1-\half \tilde{\beta}_{\rho,\epsilon}  t }\frac{\epsilon'}{\sqrt{5} \|G_0\|}
    &= \log_{\frac{1}{1-\half \tilde{\beta}_{\rho,\epsilon}  t }}\frac{\sqrt{5}\|G_0\|}{\epsilon'} \\
    &\leq \log_{\frac{1}{1-\frac{1}{4} \beta  t }}\frac{\sqrt{5}\|G_0\|}{\epsilon'} \leq \log_{\frac{4}{4-\beta t }}\frac{\sqrt{5}\omega_\rho \epsilon}{\epsilon'},
  \end{align*}
  which gives the upper bound for $T \leq T'$ stated in the theorem.

  Next, let us show that $x^T$ satisfies \eqref{eq.epsbeta} with $(\epsilon,\beta)$ replaced by $(\epsilon',\beta')$. Note that $\|G_T\| \leq \frac{\epsilon'}{\sqrt{5}}$ by the termination condition in \eqref{eq.fl.gd} and $\lambda_n(\nabla^2 F_\rho(x^T)) \geq \tilde\beta_{\rho,\epsilon}$ by \eqref{ineq.strong.convex.2}. Therefore, \cite[Corollary 2.7]{GoyeEfteBoum2024} applied at $x^T$ gives
  \bequation\label{eq.new}
  \baligned
    &\ \|c(x^T)\| \leq \|G_T\|,\ \ \|\nabla_x L(x^T,y(x^T))\| \leq 2\|G_T\|, \\ \text{and}\ 
    &\ d^T\nabla_{xx}^2 L(x^T,y(x^T))d \geq \(\tilde\beta_{\rho,\epsilon} - \(1 + \frac{2}{\sigma_{\min}}\) \kappa \|G_T\|\) \|d\|^2 \\ \text{for all}\ 
    &\ d \in \Null(\nabla c(x^T)^T).
  \ealigned
  \eequation
  From this, the termination condition in~\eqref{eq.fl.gd}, and the properties~\eqref{ineq.induct.end},
  \begin{align*}
    &\ \|\nabla L(x^{T},y(x^{T}))\| \nonumber \\
    =&\ \sqrt{\|\nabla_x L(x^{T},y(x^{T}))\|^2+\|\nabla_{\yhat} L(x^{T},y(x^{T}))\|^2} \nonumber \\
    \leq&\ \sqrt{5}\|G_T\| \leq \sqrt{5} \(\frac{\epsilon'}{\sqrt{5}}\) = \epsilon'.
  \end{align*}
  At the same time, combining \eqref{eq.beta.eps.rho}, \eqref{eq.new}, and $\|G_T\| \leq \frac{\epsilon'}{\sqrt{5}}$ one has
  \bequationNN
    \begin{aligned}
      &\ \tilde{\beta}_{\rho,\epsilon} - \(1+\frac{2}{\sigma_{\min}}\) \kappa \|G_T\| \\
      \geq &\ \tilde{\beta}_{\rho,\epsilon} - \(1+\frac{2}{\sigma_{\min}}\) \frac{\kappa}{\sqrt{5}} \epsilon' \\
      = &\ \beta - \eta_\rho\epsilon -\frac{4 M_\rho \omega_\rho}{\beta}\epsilon- \(1+\frac{2}{\sigma_{\min}}\) \frac{\kappa}{\sqrt{5}} \epsilon' \ge\half\beta.
    \end{aligned}
  \eequationNN
  which gives \eqref{eq.beta'}.  Overall, $x^T$ satisfies \eqref{eq.epsbeta} with $(\epsilon',\beta')$, as claimed.

  Lastly, let us explain that \cite[Algorithm 1]{GoyeEfteBoum2024} with the parameter choice in \eqref{goyens.alg.parameter} reduces to the gradient descent method \eqref{eq.fl.gd}. In contrast to \eqref{eq.fl.gd} itself, \cite[Algorithm 1]{GoyeEfteBoum2024} first moves along negative gradient directions until $\|\nabla F_\rho\| \leq \epsilon_1$. Then, if the curvature is sufficiently negative, it starts a second-order phase by moving along a direction computed by using second-order information. The stepsizes for both directions are determined by using a line-search strategy. From \eqref{eq.beta.eps.rho} and \eqref{ineq.strong.convex.2}, one can observe $\lambda_n(\nabla^2 F_\rho(x^{i+1})) \geq \half \beta > 0$ for all $i$, and hence the second-order phase will never be triggered. Moreover, let us show that $\lambda_1(\nabla^2 F_\rho(x)) \leq \lambda_\rho$ for all $x$ on the line segment between $x^i$ and $x^{i+1}$.  For any $s \in [0,1]$, let $x(s) = x^i + s(x^{i+1} - x^i)$. Then, by \cite[Theorem~6.6]{Stew1973}, Assumption~\ref{ass.Lagrangian}, $x^i \in \Ccal_R$ from \eqref{ineq.induct.end} along with the definition of $\lambda_{R,\rho}$ in \eqref{sigma.max.hess.F}, the fact that $\|x^{i+1} - x^i\| = t\|G_i\|$, \eqref{ineq.g0}, the fact that $\tilde\beta_{\rho,\epsilon} \leq \beta$ by \eqref{eq.beta.eps.rho}, and $\lambda_\rho \geq \beta + \lambda_{R,\rho}$, one finds that
  \begin{align*}
    \lambda_1(\nabla^2 F_\rho(x(s)))
      &\leq \lambda_1(\nabla^2 F_\rho(x^i)) + \|\nabla^2 F_\rho(x(s)) - \nabla^2 F_\rho(x^i)\| \\
      &\leq \lambda_{R,\rho} + M_\rho s \|x^{i+1} - x^i\| \\
      &\leq \lambda_{R,\rho} + M_\rho t \|G_i\| \\
      &\leq \lambda_{R,\rho} + \tilde\beta_{\rho,\epsilon} \leq \lambda_{R,\rho} + \beta \leq \lambda_\rho.
  \end{align*}
  Therefore, from this bound, Taylor's theorem, and $t$ in \eqref{eq.stepsize}, one has
  \begin{align}
    F_\rho(x^{i+1})
      &\leq F_\rho(x^i) + \nabla F_\rho(x^i)^T (x^{i+1} - x^i)  + \half \lambda_\rho  \|x^{i+1}-x^{i}\|^2 \nonumber \\
      &= F_\rho(x^i) - t\|\nabla F_\rho(x^i)\|^2  + \half \lambda_\rho  t^2\|\nabla F_\rho(x^i)\|^2  \nonumber \\
  &\le F_\rho(x^i) - \half t\|\nabla F_\rho(x^i)\|^2 \nonumber\\
  &= F_\rho(x^i) - c_1 \alpha_{01}\|\nabla F_\rho(x^i)\|^2. \label{F.taylor}
  \end{align}
  Hence, the line-search termination condition in iteration $i$ of \cite[Algorithm~\ref{alg.sub.solver}]{GoyeEfteBoum2024} is satisfied at $x^{i+1}$, so \cite[Algorithm 1]{GoyeEfteBoum2024} reduces to gradient descent with a fixed step size $t$. Moreover, the choice $\epsilon_1$ in \eqref{goyens.alg.parameter} ensures the output of \cite[Algorithm 1]{GoyeEfteBoum2024}, say $\xtilde^T$, gives $\|\nabla F_\rho(\xtilde^T)\|\le\epsilon_1=\frac{\epsilon'}{\sqrt{5}}$, which is exactly the termination condition in~\eqref{eq.fl.gd}.
  \qed
\eproof

\section{Application: Progressive Sampling}\label{sec.progressive}

The purpose of this section is to show that Theorem~\ref{theo.fl.strongly.convex} can be leveraged to prove a strong complexity guarantee in a specific setting of interest. In particular, we consider the setting of progressive sampling as proposed in \cite{CurtGuoRobi2025} to solve instances of problem~\eqref{eq.problem} when the objective and constraint functions are defined through expectations, which in turn are approximated by a large-scale sample average approximation (SAA) strategy. We show that when~\cite[Algorithm 1]{GoyeEfteBoum2024} is employed as the subproblem solver, an improved sample complexity can be obtained by a progressive sampling strategy to solve the SAA problem as compared to an approach that solves the full SAA problem directly.

The SAA problem type that we consider in this section can be written as
\bequation\label{prob.opt.N}
  \begin{aligned}
    \min_{x \in \R{n}}\ &\ f(x) \equiv \frac{1}{N} \sum_{i \in [N]} f_i(x)\\
    \st\ &\ c(x) \equiv \frac{1}{N} \sum_{i \in [N]} c_i(x) = 0,
  \end{aligned}
\eequation
where, among other properties, the objective function $f : \R{n} \to \R{}$ and constraint function $c: \R{n} \to \R{m}$ satisfy Assumption~\ref{ass.smooth}.  Generally speaking, the objective and constraint functions could be defined with respect to different numbers of terms, but to simplify our notation we assume that $f$ and $c$ are each composed of $N$ terms.  A main challenge of solving instances of problem~\eqref{prob.opt.N} is that $N$ can be quite large, meaning that each evaluation of $f$, $c$, or either of their first-order derivatives requires taking a sum over $N$ terms.

The progressive sampling strategy proposed in~\cite{CurtGuoRobi2025} has been designed for such a context.  The main hope is to exploit the fact that many of the objective and constraint function terms may be relatively similar to each other, say, in settings when the objective and constraints are averages over functions defined by independent and identically distributed draws of a random variable.  The main idea of the approach is to solve \eqref{prob.opt.N} in a progressive manner by first selecting a subset of terms (e.g., at random), solving the corresponding sampled problem to moderate accuracy, then increasing the sample size and tightening the accuracy tolerance in an iterative manner, where in each case the solve is warm-started by the approximate solution obtained for the previous sample. Denoting $\Scal_k \subseteq [N]$, a corresponding approximate objective function as $f_{\Scal_k} : \R{n} \to \R{}$, and a corresponding approximate constraint function as $c_{\Scal_k} : \R{n} \to \R{m}$, one can write a subproblem of progressive sampling as
\bequation\label{prob.opt.S}
  \baligned
    \min_{x \in \R{n}}\ f_{\Scal_k}(x)\ &\st\ c_{\Scal_k}(x) = 0, \\ \text{where}\ \ f_{\Scal_k}(x) = \frac{1}{|\Scal_k|} \sum_{i\in{\Scal_k}} f_i(x) \ \ &\text{and}\ \ c_{\Scal_k}(x) = \frac{1}{|\Scal_k|} \sum_{i\in{\Scal_k}} c_i(x).
  \ealigned
\eequation

A statement of the progressive sampling method from \cite{CurtGuoRobi2025} (assuming equal sample sizes for the objective and constraint functions) is given as Algorithm~\ref{alg.psm}.  Also, a simplified statement of \cite[Algorithm 1]{GoyeEfteBoum2024} is given as Algorithm~\ref{alg.sub.solver}.  (The step-size selection mechanisms that may be employed by the algorithm are specified in our subsequent analysis.) Analogously to Definition~\ref{def.fal}, for any $\Scal \subseteq [N]$ the algorithm defines $F_{\Scal,\rho} : \hat\Ccal_{\Scal,R} \to \R{}$ (see Assumption~\ref{ass.app.licq} below for the definition of $\hat\Ccal_{\Scal,R}$), which for all $x \in \hat\Ccal_{\Scal,R}$ is given by
\bequationNN
  F_{\Scal,\rho}(x) = f_{\Scal}(x) + c_{\Scal}(x)^Ty_{\Scal}(x) + \rho \|c_{\Scal}(x)\|^2,
\eequationNN
where analogously to \eqref{eq.y}, for any $x \in \hat\Ccal_{\Scal,R}$, the multiplier is given by
\bequation\label{eq.yS}
  y_\Scal(x) = -\nabla c_\Scal(x)^\dag \nabla f_\Scal(x).
\eequation
Let us also introduce for our analysis, as in \eqref{eq.Lagrangian}, the Lagrangian
\bequationNN
  L_\Scal(x,\yhat) = f_\Scal(x) + c_\Scal(x)^T\yhat.
\eequationNN
For these quantities, we have $F \equiv F_{[N]}$, $y \equiv y_{[N]}$, and $L \equiv L_{[N]}$.

\balgorithm
  \caption{Progressive Sampling Algorithm for Solving \eqref{prob.opt.N}}
  \label{alg.psm}
  \balgorithmic[1]
    \Require Initial sample size $p_1 \in [N]$, sample increase factor $\theta \in (1,\infty)$, initial point $x_0 \in \R{n}$, iteration limit $K = \lceil \log_{\theta} \frac{N}{p_1} \rceil + 1$, and tolerances $\{(\epsilon_k,\zeta_k)\}_{k=1}^K \subset \R{}_{>0} \times \R{}_{>0}$
    \State set $\Scal_0 \gets \emptyset$
    \For {$k \in [K]$}
      \State choose $\Scal_k \supseteq \Scal_{k-1}$ with $|\Scal_k| = p_k$
      \State choose $\rho_k \in \R{}_{>0}$
      \State using $x_{k-1}$, $\Scal_k$, $\rho_k$, and $(\epsilon_k,\zeta_k)$, call Algorithm~\ref{alg.sub.solver} to obtain $x_k$
      \State set $p_{k+1} \gets \min\{\theta p_k, N\}$
    \EndFor
    \State \Return $(x_K,y(x_K))$
  \ealgorithmic
\ealgorithm

\balgorithm
  \caption{Adapted \cite[Algorithm 1]{GoyeEfteBoum2024} for Solving \eqref{prob.opt.S}}
  \label{alg.sub.solver}
  \balgorithmic[1]
    \Require Initial point $x_{k-1} \in \R{n}$, sample set $\Scal_k \subseteq [N]$, penalty parameter $\rho_k \in \R{}_{>0}$, and tolerances $(\epsilon_k,\zeta_k)\in \R{}_{>0} \times \R{}_{>0}$ from Algorithm~\ref{alg.psm}
    \State set $i \gets 0$ and $x_{k-1}^0 \gets x_{k-1}$.
    \For {$i \in \N{}$}
      \If{$\|\nabla F_{\Scal_k,\rho_k}(x_{k-1}^i)\| > \epsilon_k$}
      \State set $x_{k-1}^{i+1} \gets x_{k-1}^{i} - t_g\nabla F_{\Scal_k,\rho_k}(x_{k-1}^{i})$ for some $t_g \in \R{}_{>0}$
      \ElsIf{$\lambda_n(\nabla^2 F_{\Scal_k,\rho_k}(x_{k-1}^i))<-\zeta_k$}
      \State set $d \in \R{n}$ with $\|d\|=1$, $d^T\nabla^2 F_{\Scal_k,\rho_k}(x_{k-1}^i)d<-\zeta_k$, and $d^T\nabla F_{\Scal_k,\rho_k}(x_{k-1}^i)\leq0$       
       \State set $x_{k-1}^{i+1} \gets x_{k-1}^{i} + t_Hd$ for some $t_H \in \R{}_{>0}$
       \Else
       \State set $x_k \gets x_{k-1}^i$
       \State \Return $x_{k}$
      \EndIf
    \EndFor
  \ealgorithmic
\ealgorithm

For our analysis of Algorithm~\ref{alg.psm}, we make the following assumptions.  (These assumptions are consistent with those in \cite{CurtGuoRobi2025}, so making them allows us to employ the results up through Theorem 3.8 in that paper.) First, extending Assumption~\ref{ass.smooth}, we assume the following.

\bassumption\label{ass.app.boundness}
  The objective function $f : \R{n} \to \R{}$ and constraint function $c : \R{n} \to \R{m}$ are twice-continuously differentiable.  Moreover, there exists $(\kappa_{\nabla f}, \kappa_{\nabla c}, \kappa_{\nabla^2 f}, \kappa_{\nabla^2 c}) \in \R{}_{>0} \times \R{}_{>0} \times \R{}_{>0} \times \R{}_{>0}$ such that, for all $x \in \R{n}$ and $j \in [m]$, one has $\| \nabla f(x) \| \leq \kappa_{\nabla f}$, $\| \nabla c(x) \| \leq \kappa_{\nabla c}$, $\| \nabla^2 f(x) \| \leq \kappa_{\nabla^2 f}$, and $\| \nabla^2 [c]_j (x) \| \leq \kappa_{\nabla^2 c}$, where recall that $[c]_j$ is the $j$th component of $c$.
\eassumption

To connect Assumption~\ref{ass.app.boundness} to properties of the approximation functions, we make the following loose assumption about the component functions. The main consequence of this assumption can be seen in \cite[Lemma~3.2]{CurtGuoRobi2025}, which shows that by increasing the sample size one obtains approximation function and derivative values that more closely approximate those of $f$ and~$c$.

\bassumption\label{ass.app.bounded.distribute}
  There exists $(\gamma_f, \gamma_c, \gamma_{\nabla f}, \gamma_{\nabla c}, \gamma_{\nabla^2 f}, \gamma_{\nabla^2 c}) \in \R{}_{>0} \times \R{}_{>0} \times \R{}_{>0} \times \R{}_{>0} \times \R{}_{>0} \times \R{}_{>0}$ such that, for all $(x,j) \in \R{n} \times [m]$, one has
  \begin{align*}
    \frac{1}{N} \sum_{i \in [N]} |f_i(x) - f(x) |^2 &\leq \gamma_f, \\
    \frac{1}{N} \sum_{i \in [N]} \| c_i(x) - c(x) \|^2 &\leq \gamma_c, \\
    \frac{1}{N} \sum_{i \in [N]} \|\nabla f_i(x) - \nabla f(x) \|^2 &\leq \gamma_{\nabla f}, \\
    \frac{1}{N} \sum_{i \in [N]} \| \nabla c_i(x)-\nabla c(x) \|^2 &\leq \gamma_{\nabla c} \|\nabla c (x) \|^2, \\
    \frac{1}{N} \sum_{i \in [N]} \| \nabla^2 f_i(x) - \nabla^2 f(x) \|^2 &\leq \gamma_{\nabla^2 f} \\
    \text{and}\ \ 
    \frac{1}{N} \sum_{i \in [N]} \| \nabla^2 [c_i]_j(x) - \nabla^2 [c]_j(x) \|^2 &\leq \gamma_{\nabla^2 c}.
  \end{align*}
\eassumption

Next, extending Assumptions~\ref{ass.compact}, \ref{ass.licq}, and \ref{ass.Lagrangian} for any sample problem encountered by the algorithm, we make the following assumptions.

\bassumption\label{ass.app.compact}
  There exists $R \in \R{}_{>0}$ such that, for any $\Scal \subseteq [N]$ with $|\Scal| \geq p_1$, the set $\Ccal_{\Scal,R} := \{x \in \R{n} : \|c_\Scal(x)\| \leq R\}$ is compact.
\eassumption

\bassumption\label{ass.app.licq}
  There exists $\delta \in \R{}_{>0}$ $($uniform over $\Scal \subseteq [N]$ with $|\Scal| \geq p_1$$)$, a bounded open set $\hat\Ccal_{\Scal,R}$ containing $\Ccal_{\Scal,R} + \{v \in \R{n} : \|v\| \leq \delta\}$, and a positive real number $\sigma_{\min} \in \R{}_{>0}$ such that, for any $\Scal \subseteq [N]$ with $|\Scal| \geq p_1$ and $x \in \hat\Ccal_{\Scal,R}$, the constraint Jacobian has $\sigma_m(\nabla c_\Scal(x)^T) \geq \sigma_{\min}$.
\eassumption

\bassumption\label{ass.app.Lagrangian}
  Along with Assumptions~\ref{ass.app.boundness}, \ref{ass.app.bounded.distribute}, \ref{ass.app.compact}, and \ref{ass.app.licq}, for any $\Scal \subseteq [N]$ with $|\Scal| \geq p_1$ the function $y_\Scal$ defined by \eqref{eq.yS} is twice-continuously differentiable over $\hat\Ccal_{\Scal,R}$, and over $\hat\Ccal_{\Scal,R}$ the functions $f_\Scal + c_\Scal^Ty_\Scal$ and $\|c_\Scal\|^2$ have Hessian functions that are Lipschitz continuous in the sense that there exists $(M_f,M_c) \in \R{}_{>0} \times \R{}_{>0}$ $($independent of~$\Scal$$)$ such that, for all $(x,\xbar) \in \hat\Ccal_{\Scal,R} \times \hat\Ccal_{\Scal,R}$, one has that
  \begin{align*}
    \|\nabla^2 (f_\Scal + c_\Scal^Ty_\Scal)(x) - \nabla^2 (f_\Scal + c_\Scal^Ty_\Scal)(\xbar)\| &\leq M_f \|x - \xbar\| \\ \text{and}\ \ 
    \|\nabla^2 (\|c_\Scal\|^2)(x) - \nabla^2 (\|c_\Scal\|^2)(\xbar)\| &\leq M_c \|x - \xbar\|.
  \end{align*}
  Under these conditions, for any $\rho \in \R{}_{>0}$, define $M_\rho := M_f + \rho M_c$.
\eassumption

Our last assumption is the one that allows us to leverage the local-linear convergence rate of Algorithm~\ref{alg.sub.solver} in the neighborhood of certain stationary points (as shown in the previous section) in order to prove an overall improved sample complexity bound for a progressive sampling strategy. It builds on the concept of a Morse function~\cite{milnor1963morse} to that of a Morse program~\cite{fujiwara1982morse}. As discussed in \cite[Remark~3.2]{CurtGuoRobi2025}, the assumption is relatively mild in the context of constrained optimization problems with isolated local minima~\cite{guillemin1975topology}.

\bassumption\label{ass.app.subproblems}
  For some pair of positive real numbers $(\alpha,\beta) \in \R{}_{>0} \times \R{}_{>0}$, problem~\eqref{prob.opt.N} is $(\alpha,\beta)$-strongly Morse in the sense that, for any $x \in \R{n}$, if the pair $(x,y(x))$ satisfies $\| \nabla L(x,y(x)) \| \leq \alpha$, then
  \bequationNN
    |d^T \nabla_{xx}^2 L(x,y(x)) d| \geq \beta \|d\|^2\ \ \text{for all}\ \ d \in \Null(\nabla c(x)^T).
  \eequationNN
\eassumption

For our analysis, we also introduce the useful quantity
\bequation\label{eq.key_value}
  \xi_\Scal := \sqrt{\frac{N(N-|\Scal|)}{|\Scal|^2}} \in [0,\sqrt{N(N-1)}] \ \ \text{for all nonempty $\Scal\subseteq [N]$.}
\eequation

The lemma below shows conditions under which we can carry over the theoretical guarantees from \cite{GoyeEfteBoum2024} for our employment of Algorithm~\ref{alg.sub.solver} within the context of Algorithm~\ref{alg.psm}.  In the lemma, we use $\tau_{\Scal_1} \in \R{}_{>0}$ to denote the positive real number ``$C$'' introduced in \cite[Corollary~2.7]{GoyeEfteBoum2024}, which depends on second-order derivatives of $\{[c_{\Scal_1}]_j\}_{j\in[m]}$ and $\{[y_{\Scal_1}]_j\}_{j\in[m]}$ over $\Ccal_{\Scal_1,R}$.  We use this notation to indicate the dependence of this quantity on the sample set $\Scal_1$.

\begin{lemma}\label{lemma.fletures.lag.result}
  Suppose that Assumptions~\ref{ass.app.boundness}, \ref{ass.app.bounded.distribute}, \ref{ass.app.compact}, \ref{ass.app.licq} and \ref{ass.app.Lagrangian} hold, a sample set $\Scal_1 \subseteq [N]$ is given with $|\Scal_1| \geq p_1$, and an initial point $x_0 \in \Ccal_{\Scal_1,R}$ is given.  Then, there exists $\hat\rho_1 \in \R{}_{>0}$ such that the requirements of \cite[Theorem 3.5]{GoyeEfteBoum2024} hold for all $\rho_1 \geq \hat\rho_1$.  Consequently, for any such $\rho_1$ there exists $(u_{\Scal_1,1},u_{\Scal_1,2}) \in \R{}_{>0} \times \R{}_{>0}$ such that for any $(\bar\epsilon,\bar\zeta) \in (0,\frac{\sqrt{5}}{2}R] \times (0,1]$ one has that Algorithm~\ref{alg.sub.solver} with inputs~$x_0$, $\Scal_1$, $\rho_1$, and $(\epsilon_1,\zeta_1) \in \R{}_{>0} \times \R{}_{>0}$ satisfying
  \bequation\label{eq.initial_max}
    \baligned
      \epsilon_1 \leq \bar\epsilon_1(\bar\epsilon,\bar\zeta) &:= \min \left\{\frac{\sqrt{5}\bar\epsilon}{5}, \frac{\bar\zeta}{4\tau_{\Scal_1}} \right\} \\
      \text{and}\ \ \zeta_1 \leq \bar\zeta_1(\bar\epsilon,\bar\zeta) &:= \frac{\bar\zeta}{2} + \tau_{\Scal_1} \bar\epsilon_1(\bar\epsilon,\bar\zeta)
    \ealigned
  \eequation
  $($using either the line-search methods in \cite[Algorithms~2--3]{GoyeEfteBoum2024} or setting the step sizes according to the lower bounds in \cite[Lemmas~3.2 and 3.4]{GoyeEfteBoum2024}$)$ terminates in a number of iterations that is at most
  \bequation\label{eq.fleture.complexity.second}
    T_{\Scal_1} = \max\{u_{\Scal_1,1} \epsilon_1^{-2}, u_{\Scal_1,2} \zeta_1^{-3}\}.
  \eequation
  Moreover, the final iterate is $(\bar\epsilon, \bar\zeta)$-stationary $($see Definition~\ref{def.approximate.stationary}$)$.
\end{lemma}
\bproof
  Consider arbitrary $\Scal_1$ and $x_0$ as stated.  Our first aim is to prove that the requirements of \cite[Theorem 3.5]{GoyeEfteBoum2024} hold for all $\rho_1$ above a threshold $\hat\rho_1 \in \R{}_{>0}$. This requires showing that Assumptions~A1--A5 in \cite{GoyeEfteBoum2024} hold. Assumption~A1 requires the existence of $\Rhat \in \R{}_{>0}$ and $\hat\sigma \in \R{}_{>0}$ such that for all $x \in \R{n}$ with $\|c_{\Scal_1}(x)\| \leq \Rhat$ one has $\sigma_{\min}(\nabla c_{\Scal_1}(x)^T) \geq \hat\sigma$.  This holds in our present setting with Assumption~\ref{ass.app.licq}.  Assumption~A2 requires the existence of $R$ as stated in Assumption~\ref{ass.app.compact}. Assumption~A3 requires the existence of $\eta_{\Scal_1} \in \R{}_{>0}$ such that
  \bequation\label{ineq.A3}
    \|c_{\Scal_1}(x + d) - c_{\Scal_1}(x) - \nabla c_{\Scal_1}(x)^Td\| \leq \eta_{\Scal_1} \|d\|^2\ \ \text{for all}\ \ (x,d) \in \Ccal_{\Scal_1,R} \times \R{n}.
  \eequation
  To show that this holds here, observe that for all $(x,j) \in \R{n} \times [m]$ it follows from Assumption~\ref{ass.app.bounded.distribute}, \cite[(3.3f)]{CurtGuoRobi2025}, and the triangle inequality that
  \begin{align}
    \|\nabla^2 [c_{\Scal_1}]_j(x)\|
      &=    \|\nabla^2 [c]_j(x) + \nabla^2 [c_{\Scal_1}]_j(x) - \nabla^2 [c]_j(x)\| \nonumber \\
      &\leq \|\nabla^2 [c]_j(x)\| + \|\nabla^2 [c_{\Scal_1}]_j(x) - \nabla^2 [c]_j(x)\| \nonumber \\
      &\leq \|\nabla^2 [c]_j(x)\| + \xi_{\Scal_1} \sqrt{\gamma_{\nabla^2 c}} \nonumber \\
      &\leq \kappa_{\nabla^2 c} + \sqrt{\frac{N(N-p_1)\gamma_{\nabla^2c}}{p_1^2}} =: \kappa_{p_1}. \label{ineq.used}
  \end{align}
  Now consider arbitrary $(x,d) \in \Ccal_{\Scal_1,R} \times \R{n}$ and observe that, by Taylor's theorem and \eqref{ineq.used}, it follows that for all $j \in [m]$ there exists $\xtilde_j \in \R{n}$ with
  \begin{align}
    &\ |[c_{\Scal_1}]_j(x + d) - [c_{\Scal_1}]_j(x) - \nabla [c_{\Scal_1}]_j(x)^Td| \nonumber \\
    =&\ \half |d^T \nabla^2 [c_{\Scal_1}]_j(\xtilde_j)d| \leq \half \|\nabla^2 [c_{\Scal_1}]_j(\xtilde_j)\| \|d\|^2 \leq \half \kappa_{p_1} \|d\|^2. \label{ineq.taylor}
  \end{align}
  Consequently, one finds that
  \bequationNN
    \|c_{\Scal_1}(x + d) - c_{\Scal_1}(x) - \nabla c_{\Scal_1}(x)^Td\| \leq \sqrt{\sum_{j \in [m]} \frac{1}{4} \kappa_{p_1}^2 \|d\|^4} = \half \sqrt{m} \kappa_{p_1} \|d\|^2,
  \eequationNN
  from which it follows that \eqref{ineq.A3} holds for $\eta_{\Scal_1} = \half \sqrt{m} \kappa_{p_1}$.  Assumption~A4 requires that $x_0 \in \Ccal_{\Scal_1,R}$, as is stated for this lemma.  Finally, Assumption~A5 requires that the penalty parameter is chosen greater than
  \bequationNN
    \max_{x\in\Ccal_{\Scal_1,R}} \max \left\{ \frac{\sigma_1(\nabla c_{\Scal_1}(x)) \sigma_1(\nabla y_{\Scal_1}(x))}{2\sigma_m(\nabla c_{\Scal_1}(x))^2}, \frac{\sigma_1(\nabla y_{\Scal_1}(x))}{\sigma_m(\nabla c_{\Scal_1}(x))}, \frac{1}{\sigma_m(\nabla c_{\Scal_1}(x))} \right\}.
  \eequationNN
  Since $\Ccal_{\Scal_1,R}$ is compact, it follows under Assumption~\ref{ass.app.licq} and the extreme value theorem that there exists $\hat\rho_1$ as stated in the lemma.  All together, one can conclude that the requirements of \cite[Theorem 3.5]{GoyeEfteBoum2024} hold in our present setting.

  Next we prove the claimed worst-case complexity bound and the claimed $(\bar\epsilon,\bar\zeta)$-stationarity of the final iterate produced by Algorithm~\ref{alg.sub.solver}.  Toward these ends, one first has from \cite[Theorem 3.5]{GoyeEfteBoum2024} that there exists $(\ubar_{\Scal_1,1}, \ubar_{\Scal_1,2}, \ubar_{\Scal_1,3}) \in \R{}_{>0} \times \R{}_{>0} \times \R{}_{>0}$ such that in a number of iterations that is at most
  \bequation\label{eq.TS}
    \overline{T}_{\Scal_1} := \max \left\{ \ubar_{\Scal_1,1} \epsilon_1^{-2}, \ubar_{\Scal_1,2} \zeta_1^{-1}, \ubar_{\Scal_1,3} \zeta_1^{-3} \right\}
  \eequation
  the algorithm produces a point $\xbar \in \R{n}$ satisfying
  \bequation\label{ineq.fletcher.termin.second}
    \|\nabla F_{\Scal_1}(\xbar)\| \leq \epsilon_1\ \ \text{and}\ \ d^T\nabla^2 F_{\Scal_1}(\xbar) d \geq -\zeta_1 \|d\|^2\ \ \text{for all}\ \ d \in \R{n},
  \eequation
  meaning that Algorithm~\ref{alg.sub.solver} terminates with $\xbar$ as the final iterate.  (Note that to apply~\cite[Theorem 3.5]{GoyeEfteBoum2024} it has been observed that $\epsilon_1 \leq \frac{\sqrt{5}}{5} \bar\epsilon \leq \half R$.)  To show that this is achieved in a number of iterations that is at most $T_{\Scal_1}$ in~\eqref{eq.fleture.complexity.second}, one only needs to observe that $\bar\zeta \in (0,1]$ implies
  \bequationNN
    \zeta_1 \leq \frac{\bar\zeta}{2} + \tau_{\Scal_1} \(\frac{\bar\zeta}{4 \tau_{\Scal_1}}\) \leq \bar\zeta \leq 1 \implies \zeta_1^{-1} \leq \zeta_1^{-3},
  \eequationNN
  so from \eqref{eq.TS} the bound in \eqref{eq.fleture.complexity.second} holds with $u_{\Scal_1,1} := \ubar_{\Scal_1,1}$ and $u_{\Scal_1,2} := \max\{\ubar_{\Scal_1,2},\ubar_{\Scal_1,3}\}$.  Finally, \cite[Theorem 3.5]{GoyeEfteBoum2024} also shows that $\xbar$ satisfies
  \begin{align*}
    \|\nabla_x L(\xbar,y(\xbar))\| &\leq 2\epsilon_1,\ \ \|\nabla_{\yhat} L(\xbar,y(\xbar))\| \leq \epsilon_1, \\ \text{and}\ \ 
    d^T\nabla_{xx}^2 L(\xbar,y(\xbar))d &\geq -(\zeta_1 + \tau_{\Scal_1}\epsilon_1) \|d\|^2\ \ \text{for all}\ \ d \in \Null(\nabla c_{\Scal_1}(\xbar)^T).
  \end{align*}
  Observing that this means $\|\nabla L(\xbar,y(\xbar))\| \leq \sqrt{5} \epsilon_1 \leq \bar\epsilon$, and observing that $\zeta_1 + \tau_{\Scal_1} \epsilon_1 \leq \frac{\bar\zeta}{2} + 2\tau_{\Scal_1} (\frac{\bar\zeta}{4\tau_{\Scal_1}}) \leq \bar\zeta$, it follows that $\xbar$ is $(\bar\epsilon,\bar\zeta)$-stationary.
  \qed
\eproof

We now state a useful result, namely, \cite[Lemma 3.6]{CurtGuoRobi2025}, which shows that the strong Morse property extends from the full-sample problem to any sampled problem as long as the sample size is sufficiently large.

\begin{lemma}\label{lem.strongMorseextend}
  (\cite[Lemma 3.6]{CurtGuoRobi2025}) Suppose that Assumptions~\ref{ass.app.boundness}, \ref{ass.app.bounded.distribute}, \ref{ass.app.compact}, \ref{ass.app.licq}, \ref{ass.app.Lagrangian}, and~\ref{ass.app.subproblems} hold, and define the tuple $(\kappa_1, \kappa_2, \kappa_3) \in \R{}_{>0} \times \R{}_{>0} \times \R{}_{>0}$ by
  \begin{align*}
    &\ \kappa_1 := \frac{\sqrt{\gamma_{\nabla c}} \kappa_{\nabla c}}{\sigma_{\min}},\ \ \kappa_2 := \sqrt{(2 \sqrt{\kappa_{\nabla f}} + \kappa_1 \kappa_{\nabla f})^2 + \gamma_c},\ \ \text{and}\ \ \\
    &\ \kappa_3 := \(3\kappa_{\nabla^2 f}\kappa_1 + \sqrt{\kappa_{\nabla^2 f}} + \frac{\sqrt{m} (3 \kappa_{\nabla^2 c} \kappa_1 + \sqrt{\gamma_{\nabla^2 c}}) (5 \kappa_{\nabla f} \kappa_1 + \sqrt{\gamma_{\nabla f}})}{2 \sigma_{\min} \kappa_1}\).
  \end{align*}
  Then, with $(\alpha,\beta) \in \R{}_{>0} \times \R{}_{>0}$ from Assumption~\ref{ass.app.subproblems}, if $\Scal_{k+1} \subseteq [N]$ has
  \bequation\label{eq.lemma2.S}
    \xi_{\Scal_{k+1}} \leq \min \left\{ \frac{1}{3\kappa_1}, \frac{\alpha}{2\kappa_2}, \frac{7\beta}{18\kappa_3} \right\},
  \eequation
  then \eqref{prob.opt.S} with $k$ replaced by $k+1$ is $(\alpha_{\Scal_{k+1}},\beta_{\Scal_{k+1}})$-strongly Morse with
  \begin{align*}
    \alpha_{\Scal_{k+1}} &:= \alpha - \kappa_2 \xi_{\Scal_{k+1}} \geq \half \alpha \\ \text{and}\ \ 
    \beta_{\Scal_{k+1}} &:= \(1-\frac{1}{3} \kappa_1 \xi_{\Scal_{k+1}}\) \beta - \kappa_3 \xi_{\Scal_{k+1}}\geq \half \beta.
  \end{align*}
\end{lemma}

Our next lemma shows that if the achieved termination tolerances for a sample problem are set in an appropriate manner, then the initial point for the solve of the next sample problem (if $K > 1$) satisfies certain critical bounds. Here, we first see the significance of Assumption~\ref{ass.app.subproblems}.

\begin{lemma}\label{lem.start_for_next}
  Suppose that Assumptions~\ref{ass.app.boundness}, \ref{ass.app.bounded.distribute}, \ref{ass.app.compact}, \ref{ass.app.licq}, \ref{ass.app.Lagrangian}, and \ref{ass.app.subproblems} hold and $K > 1$.  Define $(\kappa_1,\kappa_2,\kappa_3) \in \R{}_{>0} \times \R{}_{>0} \times \R{}_{>0}$ as in Lemma~\ref{lem.strongMorseextend} and define the pair of positive real numbers $(\tau_1,\tau_2) \in \R{}_{>0} \times \R{}_{>0}$ by
  \bequation\label{eq.def.tau2}
    \tau_1 := \kappa_2\ \ \text{and}\ \ \tau_2 := \( \kappa_{\nabla^2 f} + \frac{\sqrt{m} \kappa_{\nabla f} \kappa_{\nabla^2 c}}{\sigma_{\min}}\)\kappa_1.
  \eequation
  Suppose that $p_1 \in [N]$ is sufficiently large such that, for all $k \in [K]$,
  \bequation\label{eq.theorem2.S}
    \xi_{\Scal_k} \leq \xi_{\Scal_1} \leq \min \left\{ \frac{1}{3\kappa_1}, \frac{\alpha}{4\kappa_2}, \frac{2\beta}{9\kappa_3} \right\},
  \eequation
  and that, for all $k \in [K-1]$, a pair $(\hat\epsilon_k,\hat\zeta_k) \in \R{}_{>0} \times \R{}_{>0}$ satisfies
  \bequation\label{eq.tolorence}
    \hat\epsilon_k \leq \tau_1 \xi_{\Scal_1}\ \ \text{and}\ \ \hat\zeta_k \leq \tau_2 \xi_{\Scal_k}.
  \eequation
  Then, for all $k \in [K-1]$, if $x_k \in \R{n}$ is $(\hat\epsilon_k,\hat\zeta_k)$-stationary with multiplier $y_{\Scal_k}(x_k)$ respect to~\eqref{prob.opt.S}, then for~\eqref{prob.opt.S} with $k$ replaced by $k+1$ one has
  \bequation\label{ineq.init.next}
    \begin{aligned}
      \|\nabla L_{\Scal_{k+1}}(x_k, y_{\Scal_{k+1}}(x_k))\| \leq&\ 3 \tau_1 \xi_{\Scal_1} \leq \alpha_{\Scal_{k+1}}\ \ \text{and} \\
      d^T \nabla^2_{xx} L_{\Scal_{k+1}}(x_k, y_{\Scal_{k+1}}(x_k)) d \geq&\  \beta_{\Scal_{k+1}} \|d\|^2\ \text{for all}\ d \in \Null(\nabla c_{\Scal_{k+1}}(x_k)^T),
    \end{aligned}
  \eequation
  where one defines $(\alpha_{\Scal_{k+1}},\beta_{\Scal_{k+1}})$ as in Lemma~\ref{lem.strongMorseextend}.
\end{lemma}
\bproof
  Suppose that $K > 1$ and $x_k \in \R{n}$ is $(\hat\epsilon_k,\hat\zeta_k)$-stationary with respect to~\eqref{prob.opt.S} for $\Scal = \Scal_k$ with $(\hat\epsilon_k,\hat\zeta_k)$ satisfying~\eqref{eq.tolorence}. Our aim is to prove \eqref{ineq.init.next}. First let us bound the norm of the difference between the gradients of $L_{\Scal_k}$ and $L_{\Scal_{k+1}}$ with respect to the point~$x_k$.  Toward this end, using~\eqref{eq.yS}, define $y_k := y_{\Scal_k}(x_k)$ and $y_{k+1} := y_{\Scal_{k+1}}(x_k)$.  Then, by the triangle inequality,
  \begin{align}
    &\ \|\nabla_x L_{\Scal_{k+1}}(x_k,y_{k+1}) - \nabla_x L_{\Scal_k}(x_k,y_k)\| \nonumber \\
    \leq&\ \|\nabla_x L_{\Scal_{k+1}}(x_k,y_{k+1}) - \nabla_x L(x_k,y_k)\| \nonumber \\
    &\ + \|\nabla_x L(x_k,y_k)  - \nabla_x L_{\Scal_k}(x_k,y_k)\| .\label{ineq.L.grad.x.diff.2.2}
  \end{align}
  Then, by similar arguments as led to \cite[(29)]{CurtGuoRobi2025}, one finds
  \begin{align*}
    \|\nabla_x L_{\Scal_k}(x_k,y_k) - \nabla_x L(x_k,y_k)\|
    &\leq \xi_{\Scal_k} (2 \sqrt{\kappa_{\nabla f}} + \kappa_1 \kappa_{\nabla f})\ \ \text{and} \\
    \|\nabla_x L_{\Scal_{k+1}}(x_k,y_{k+1}) - \nabla_x L(x_k,y_k)\| &\leq \xi_{\Scal_{k+1}} (2 \sqrt{\kappa_{\nabla f}} + \kappa_1 \kappa_{\nabla f}) \\
    &\leq \xi_{\Scal_k} (2 \sqrt{\kappa_{\nabla f}} + \kappa_1 \kappa_{\nabla f}),
  \end{align*}
  where the last inequality follows from \eqref{eq.key_value} and $|\Scal_{k+1}| > |\Scal_k|$.  Combining these bounds with \eqref{ineq.L.grad.x.diff.2.2}, one obtains that
  \bequation\label{ineq.theorem2.Lx}
    \begin{aligned}
      \|\nabla_x L_{\Scal_{k+1}}(x_k,y_{k+1}) - \nabla_x L_{\Scal_k}(x_k,y_k)\| 
      \leq  2\xi_{\Scal_k} (2 \sqrt{\kappa_{\nabla f}} + \kappa_1 \kappa_{\nabla f}). 
    \end{aligned}
  \eequation
  On the other hand, from the triangle inequality, \cite[(10b)]{CurtGuoRobi2025}, and $\xi_{\Scal_{k+1}} \leq \xi_{\Scal_k}$,
  \begin{align}
    &\ \|\nabla_{\yhat} L_{\Scal_{k+1}}(x_k,y_{k+1}) - \nabla_{\yhat} L_{\Scal_k}(x_k,y_k)\| \nonumber \\
    =&\ \|c_{\Scal_{k+1}}(x_k) - c_{\Scal_k}(x_k)\| \nonumber \\
    \leq&\ \|c_{\Scal_{k+1}}(x_k) - c(x_k)\| + \|c(x_k) - c_{\Scal_k}(x_k)\| \nonumber \\
    \leq&\ \sqrt{\gamma_c} \xi_{\Scal_{k+1}} + \sqrt{\gamma_c} \xi_{\Scal_k} \leq 2 \sqrt{\gamma_c} \xi_{\Scal_k}. \label{ineq.theorem2.Ly}
  \end{align}
  Combining \eqref{ineq.theorem2.Lx} and \eqref{ineq.theorem2.Ly}, one finds that
  \bequation\label{ineq.LklK1}
  \begin{aligned}
     &\ \|\nabla L_{\Scal_{k+1}}(x_k,y_{k+1}) - \nabla L_{\Scal_k}(x_k,y_k)\|^2 \\
     \leq &\ 4 (\gamma_c + (2 \sqrt{\kappa_{\nabla f}} + \kappa_1 \kappa_{\nabla f})^2) \xi_{\Scal_k}^2 = 4 \kappa_2^2 \xi_{\Scal_k}^2,
  \end{aligned}
  \eequation
  which further gives under the conditions of the lemma that
  \begin{align}
    &\ \|\nabla L_{\Scal_{k+1}}(x_k,y_{k+1})\| \nonumber \\
    \leq&\ \|\nabla L_{\Scal_{k+1}}(x_k,y_{k+1}) - \nabla L_{\Scal_k}(x_k,y_k)\| + \|\nabla L_{\Scal_k}(x_k,y_k)\| \nonumber \\
    \leq&\ 2 \kappa_2 \xi_{\Scal_k} + \hat\epsilon_k \leq 2\tau_1 \xi_{\Scal_1} + \tau_1 \xi_{\Scal_1} \leq 3 \tau_1 \xi_{\Scal_1} = 3\kappa_2\xi_{\Scal_1} \leq \frac{3}{4} \alpha. \label{ineq.3/4alpha}
  \end{align}
  At the same time, by $\xi_{\Scal_{k+1}} \leq \xi_{\Scal_k}$ and $\kappa_2 \xi_{\Scal_k} \leq \frac{1}{4} \alpha$, one has that
  \bequation\label{ineq.alpha.k.k+1}
    \alpha_{\Scal_{k+1}} = \alpha - \kappa_2 \xi_{\Scal_{k+1}} \geq \alpha - \kappa_2 \xi_{\Scal_k} \geq \alpha - \frac{1}{4}\alpha = \frac{3}{4}\alpha.
  \eequation
  Combining \eqref{ineq.3/4alpha} and \eqref{ineq.alpha.k.k+1}, one obtains the first conclusion in \eqref{ineq.init.next} that
  \bequation\label{ineq.S.k+1.grad.cond}
    \|\nabla L_{\Scal_{k+1}}(x_k, y_{k+1})\| \leq 3 \tau_1 \xi_{\Scal_1} \leq \alpha_{\Scal_{k+1}}.
  \eequation

  Our next goal is to prove the second inequality in \eqref{ineq.init.next}.  Toward this end, first note that Lemma~\ref{lem.strongMorseextend} applies here since the right-hand side of \eqref{eq.theorem2.S} is less than or equal to the right-hand side of \eqref{eq.lemma2.S}.  Hence, by Lemma~\ref{lem.strongMorseextend}, it follows that~\eqref{prob.opt.S} with $k$ replaced by $k+1$ is $(\alpha_{\Scal_{k+1}},\beta_{\Scal_{k+1}})$-strongly Morse.  Combined with~\eqref{ineq.S.k+1.grad.cond}, this means that one has for all $d_{\Scal_{k+1}} \in \Null(\nabla c_{\Scal_{k+1}}(x_k)^T)$ that
  \bequation\label{eq.lemma36}
    |d_{\Scal_{k+1}}^T \nabla^2_{xx} L_{\Scal_{k+1}}(x_k,y_{k+1}) d_{\Scal_{k+1}}| \geq \beta_{\Scal_{k+1}} \|d_{\Scal_{k+1}}\|^2.
  \eequation
  To prove the second inequality in \eqref{ineq.init.next}, it is necessary to show that this inequality holds without the absolute value on the left-hand side.  Toward showing this, let us define $\dbar_{\Scal_{k+1}} := \Ncal(\nabla c_{\Scal_k}(x_k)^T) d_{\Scal_{k+1}}$.  Then, one finds
  \begin{align}
    &\ d_{\Scal_{k+1}}^T \nabla^2_{xx} L_{\Scal_{k+1}}(x_k,y_{k+1}) d_{\Scal_{k+1}} \nonumber \\
    =&\ \dbar_{\Scal_{k+1}}^T \nabla^2_{xx} L_{\Scal_k}(x_k,y_k) \dbar_{\Scal_{k+1}} \nonumber \\
    &\ + (d_{\Scal_{k+1}}^T \nabla^2_{xx} L_{\Scal_{k+1}}(x_k,y_{k+1}) d_{\Scal_{k+1}} - \dbar_{\Scal_{k+1}}^T \nabla^2_{xx} L_{\Scal_k}(x_k,y_k) \dbar_{\Scal_{k+1}}) \nonumber \\
    \geq&\ \dbar_{\Scal_{k+1}}^T \nabla^2_{xx} L_{\Scal_k}(x_k,y_k) \dbar_{\Scal_{k+1}} \nonumber \\
    &\ -|d_{\Scal_{k+1}}^T \nabla^2_{xx} L_{\Scal_{k+1}}(x_k,y_{k+1}) d_{\Scal_{k+1}} - \dbar_{\Scal_{k+1}}^T \nabla^2_{xx} L_{\Scal_k}(x_k,y_k) \dbar_{\Scal_{k+1}}|. \label{ineq.dld}
  \end{align}
  Since $\Ncal(\nabla c_{\Scal_k}(x_k)^T)$ is a projection matrix, $\|\Ncal(\nabla c_{\Scal_k}(x_k)^T)\| \leq 1$ and
  \bequationNN
  \begin{aligned}
      \|\dbar_{\Scal_{k+1}}\|^2 
      &= \|\Ncal(\nabla c_{\Scal_k}(x_k)^T)  d_{\Scal_{k+1}}\|^2 \\
      &\leq \|\Ncal(\nabla c_{\Scal_k}(x_k)^T)\|^2 \|d_{\Scal_{k+1}}\|^2 \leq \|d_{\Scal_{k+1}}\|^2.
  \end{aligned}
  \eequationNN
  Along with the fact that $x_k$ is $(\hat\epsilon_k,\hat\zeta_k)$-stationary, this means that the first-term on the right-hand side of \eqref{ineq.dld} satisfies the inequalities
  \bequation\label{ineq.varepson.cond}
    \dbar_{\Scal_{k+1}}^T \nabla^2_{xx} L_{\Scal_k}(x_k,y_k) \dbar_{\Scal_{k+1}} \geq -\hat\zeta_k \|\dbar_{\Scal_{k+1}}\|^2 \geq -\hat\zeta_k \|d_{\Scal_{k+1}}\|^2.
  \eequation
  On the other hand, with respect to the second term on the right-hand side of \eqref{ineq.dld}, observe that with the triangle inequality one finds
  \begin{align}
    &\ |d_{\Scal_{k+1}}^T \nabla^2_{xx} L_{\Scal_{k+1}}(x_k,y_{k+1}) d_{\Scal_{k+1}} - \dbar_{\Scal_{k+1}}^T \nabla^2_{xx} L_{\Scal_k}(x_k,y_k) \dbar_{\Scal_{k+1}}| \nonumber \\
    \leq&\ |d_{\Scal_{k+1}}^T \nabla^2_{xx} L_{\Scal_{k+1}}(x_k,y_{k+1}) d_{\Scal_{k+1}} - d_{\Scal_{k+1}}^T \nabla^2_{xx} L(x_k,y_k) d_{\Scal_{k+1}} | \nonumber \\
    &\ + |d_{\Scal_{k+1}}^T \nabla^2_{xx} L(x_k,y_k) d_{\Scal_{k+1}} - \dbar_{\Scal_{k+1}}^T \nabla^2_{xx} L(x_k,y_k) \dbar_{\Scal_{k+1}}| \nonumber \\
    &\ + |\dbar_{\Scal_{k+1}}^T \nabla^2_{xx} L(x_k,y_k) \dbar_{\Scal_{k+1}} - \dbar_{\Scal_{k+1}}^T \nabla^2_{xx} L_{\Scal_k}(x_k,y_k) \dbar_{\Scal_{k+1}}|. \label{ineq.expand.3.24.1}
  \end{align}
  Next, one can observe that the first and third terms on the right-hand side of~\eqref{ineq.expand.3.24.1} can be bounded as in \cite[(41)]{CurtGuoRobi2025}  with respect to $\Scal = \Scal_k$ and $\Scal = \Scal_{k+1}$, respectively.  To see this, note that a bound of the form in \cite[(36)]{CurtGuoRobi2025}  holds since \cite[(36)]{CurtGuoRobi2025} followed using general matrix-norm inequalities.  Furthermore, a bound of the form in \cite[(37)]{CurtGuoRobi2025} that relies on \cite[(36)]{CurtGuoRobi2025} and \cite[(10f)]{CurtGuoRobi2025} holds, and bounds of the form in \cite[(38)--(39)]{CurtGuoRobi2025} hold since these rely on $\xi_{\Scal_k} \leq \frac{1}{3\kappa_1}$ and $\xi_{\Scal_{k+1}} \leq \frac{1}{3\kappa_1}$, which hold here.  Thus, since \cite[(41)]{CurtGuoRobi2025} follows from a combination of \cite[(36)--(39)]{CurtGuoRobi2025}, one obtains that with
  \bequationNN
    \lambda := \sqrt{\kappa_{\nabla^2 f}} + \frac{3 \sqrt{m} \kappa_{\nabla^2 c} \kappa_1 (3\kappa_{\nabla f} \kappa_1 + \sqrt{\gamma_{\nabla f}}) + \sqrt{m} \sqrt{\gamma_{\nabla^2 c}} (5 \kappa_{\nabla f} \kappa_1 + \sqrt{\gamma_{\nabla f}})}{2 \sigma_{\min} \kappa_1}, 
  \eequationNN
  the first and third terms on the right-hand side of~\eqref{ineq.expand.3.24.1} respectively satisfy
  \begin{align}
    &\ |d_{\Scal_{k+1}}^T \nabla^2_{xx} L_{\Scal_{k+1}}(x_k,y_{k+1}) d_{\Scal_{k+1}} - d_{\Scal_{k+1}}^T \nabla^2_{xx} L(x_k,y_k) d_{\Scal_{k+1}}| \nonumber \\
    \leq&\ \lambda \xi_{\Scal_{k+1}} \|d_{\Scal_{k+1}}\|^2 \leq \lambda \xi_{\Scal_k} \|d_{\Scal_{k+1}}\|^2 \label{ineq.3.24.1} \\
    \text{and}\ \ 
    &\ |\dbar_{\Scal_{k+1}}^T \nabla^2_{xx} L(x_k,y_k) \dbar_{\Scal_{k+1}} - \dbar_{\Scal_{k+1}}^T \nabla^2_{xx} L_{\Scal_k}(x_k,y_k) \dbar_{\Scal_{k+1}}| \nonumber \\
    \leq&\ \lambda \xi_{\Scal_k} \|\dbar_{\Scal_{k+1}}\|^2 \leq \lambda \xi_{\Scal_k} \|d_{\Scal_{k+1}}\|^2. \label{ineq.3.24.2}
  \end{align}
  As for the second term on the right of \eqref{ineq.expand.3.24.1}, submultiplicity of the matrix 2-norm, the triangle inequality, $\|\dbar_{\Scal_{k+1}}\|\le\|d_{\Scal_{k+1}}\|$, and \cite[(8c)]{CurtGuoRobi2025} give
  \begin{align}
    &\ |d_{\Scal_{k+1}}^T \nabla^2_{xx} L(x_k,y_k) d_{\Scal_{k+1}} - \dbar_{\Scal_{k+1}}^T \nabla^2_{xx} L(x_k,y_k) \dbar_{\Scal_{k+1}}| \nonumber \\
    =&\ |(d_{\Scal_{k+1}} - \dbar_{\Scal_{k+1}})^T \nabla^2_{xx} L(x_k,y_k) (d_{\Scal_{k+1}} + \dbar_{\Scal_{k+1}})| \nonumber \\
    \leq&\ \|\nabla^2_{xx} L(x_k,y_k)\| \|d_{\Scal_{k+1}} - \dbar_{\Scal_{k+1}}\|\|d_{\Scal_{k+1}} + \dbar_{\Scal_{k+1}}\| \nonumber \\
    \leq&\ 2\(\kappa_{\nabla^2 f} + \frac{\sqrt{m} \kappa_{\nabla f} \kappa_{\nabla^2 c}}{\sigma_{\min}}\)\|d_{\Scal_{k+1}}\| \|d_{\Scal_{k+1}} - \dbar_{\Scal_{k+1}}\| \nonumber\\
    =&\ 2\frac{\tau_2}{\kappa_1}\|d_{\Scal_{k+1}}\| \|d_{\Scal_{k+1}} - \dbar_{\Scal_{k+1}}\|. \label{ineq.theorem2.iv2}
  \end{align}
  By Assumption~\ref{ass.app.licq} the matrices $\nabla c(x_k)$ and $\nabla c_{\Scal_k}(x_k)$ have full row rank, i.e., the same rank, which by \cite[Theorems 2.3--2.4]{Stew1977} means that
  \bequationNN
    \|(\Rcal(\nabla c(x_k)) - \Rcal(\nabla c_{\Scal_k}(x_k)))\| \leq \kappa_1 \xi_{\Scal_k}.
  \eequationNN
  Hence, with respect to $\|d_{\Scal_{k+1}}- \dbar_{\Scal_{k+1}}\|$, one finds from the triangle inequality, submultiplicity of the matrix 2-norm, and \cite[Lemma 3.5]{CurtGuoRobi2025} that
  \begin{align}\nonumber
    &\ \|d_{\Scal_{k+1}} - \dbar_{\Scal_{k+1}}\| \\
    =&\ \|\Rcal(\nabla c_{\Scal_k}(x_k)) d_{\Scal_{k+1}}\| \nonumber \\
    \leq&\ \|\Rcal(\nabla c(x_k)) d_{\Scal_{k+1}}\| + \|(\Rcal(\nabla c(x_k)) - \Rcal(\nabla c_{\Scal_k}(x_k)))\| \|d_{\Scal_{k+1}}\|\nonumber \\
    \leq&\ \kappa_1 \xi_{\Scal_{k+1}} \|d_{\Scal_{k+1}}\| + \kappa_1 \xi_{\Scal_k} \|d_{\Scal_{k+1}}\| \leq 2 \kappa_1 \xi_{\Scal_k} \|d_{\Scal_{k+1}}\|.\label{ineq.theorem2.d.dbar}
  \end{align}
  Combining \eqref{ineq.dld}--\eqref{ineq.theorem2.d.dbar} with \cite[(25)]{CurtGuoRobi2025}, \eqref{eq.def.tau2}, and $\hat\zeta_k \leq \tau_2 \xi_{\Scal_k}$, one finds
  \begin{align*}
    &\ d_{\Scal_{k+1}}^T \nabla^2_{xx} L_{\Scal_{k+1}}(x_k, y_{k+1}) d_{\Scal_{k+1}} \nonumber \\
    \geq&\ -\hat\zeta_k \|d_{\Scal_{k+1}}\|^2 - 2 \lambda \xi_{\Scal_k} \|d_{\Scal_{k+1}}\|^2 - 4\tau_2 \xi_{\Scal_k} \|d_{\Scal_{k+1}}\|^2 \\
   =&\ \(-5\tau_2  - 2 \lambda \)\xi_{\Scal_k} \|d_{\Scal_{k+1}}\|^2.
  \end{align*}
  On the other hand, from the definitions of $\tau_2$, $\lambda$, and $\kappa_3$ one finds
  \begin{align*}
    &\ \(-5\tau_2  - 2 \lambda \) \\
    \geq&\ - 2\(3\kappa_{\nabla^2 f}\kappa_1 + \sqrt{\kappa_{\nabla^2 f}} + \frac{\sqrt{m} (3 \kappa_{\nabla^2 c} \kappa_1 + \sqrt{\gamma_{\nabla^2 c}}) (5 \kappa_{\nabla f} \kappa_1 + \sqrt{\gamma_{\nabla f}})}{2 \sigma_{\min} \kappa_1}\)\\
    =&\ - 2 \kappa_3.
  \end{align*}
  Hence, since \eqref{eq.theorem2.S} requires $\xi_{\Scal_k} \leq \frac{2\beta}{9\kappa_3}$, it follows that
  \bequation\label{ineq.var.k+1}
    d_{\Scal_{k+1}}^T \nabla^2_{xx} L_{\Scal_{k+1}}(x_k,y_{k+1}) d_{\Scal_{k+1}} \geq -2 \kappa_3 \xi_{\Scal_k} \|d_{\Scal_{k+1}}\|^2 \geq -\frac{4}{9} \beta \|d_{\Scal_{k+1}}\|^2.
  \eequation
  On the other hand, recall $\beta_{\Scal_{k+1}} = \(1-\frac{1}{3} \kappa_1 \xi_{\Scal_{k+1}}\) \beta - \kappa_3 \xi_{\Scal_{k+1}}$, which along with $\xi_{\Scal_{k+1}} \leq \xi_{\Scal_k}$ and \eqref{eq.theorem2.S} (specifically, with $\xi_{\Scal_k} \leq \min\{\frac{1}{3\kappa_1}, \frac{2\beta}{9\kappa_3}\})$ implies
  \bequationNN
    \beta_{\Scal_{k+1}} \geq \(1-\frac{1}{3} \kappa_1 \xi_{\Scal_k}\) \beta - \kappa_3 \xi_{\Scal_k} \geq \beta - \frac19 \beta - \frac29 \beta = \frac{2}{3} \beta.
  \eequationNN
  Since~\eqref{prob.opt.S} with $k$ replaced by $k+1$ is $(\alpha_{k+1},\beta_{k+1})$-strongly Morse (which follows as a consequence of Lemma~\ref{lem.strongMorseextend}), one has that
  \bequationNN
    |d_{\Scal_{k+1}}^T \nabla^2_{xx} L_{\Scal_{k+1}}(x_k,y_{k+1}) d_{\Scal_{k+1}}| \geq \beta_{\Scal_{k+1}} \|d_{\Scal_{k+1}}\|^2 \geq \frac{2}{3}\beta \|d_{\Scal_{k+1}}\|^2.
  \eequationNN
  That said, since \eqref{ineq.var.k+1} holds and $-\frac49 > -\frac23$, it must hold that
  \bequationNN
    d_{\Scal_{k+1}}^T \nabla^2_{xx} L_{\Scal_{k+1}}(x_k,y_{k+1}) d_{\Scal_{k+1}} \geq \beta_{\Scal_{k+1}} \|d_{\Scal_{k+1}}\|^2,
  \eequationNN
  which completes the proof.
  \qed
\eproof

We are almost prepared to prove our main theorem about progressive sampling.  We remark upfront that, in order for our theoretical guarantee to hold for $p_1 < N$, the algorithmic parameters need to be chosen within prescribed intervals.  For notational simplicity, we shall define these intervals with respect to a single prescribed parameter $\mu \in (0,1)$. Choosing the parameters in such intervals ensures that the computational effort is balanced appropriately for all $k \in [K]$, so that, for example, the tolerance for solving one sampled problem is not arbitrarily tighter than the tolerances for solving the other sampled problems in the sequence. Such requirements can be seen to be natural for worst-case complexity guarantees of the type proved in the theorem. We emphasize, however, that for good practical performance, we do not expect these parameter requirements to be absolutely necessary. Indeed, the experiments in \cite{CurtGuoRobi2025} show good performance without such stringent parameter choices.

To define our parameter intervals, we refer to the user-defined $\mu \in (0,1)$ along with $(\alpha,\beta)$ from Assumption~\ref{ass.app.subproblems}. First, for each $\Scal \subseteq [N]$ with $|\Scal| \geq p_1$, let $\kappa_\Scal$ be defined in the same manner as $\kappa$ is defined in Lemma~\ref{lemma.F.psd}, in this case with respect to problem~\eqref{prob.opt.S} with $\Scal$ in place of $\Scal_k$, then let
\bequationNN
  \bar\kappa := \max\{ \kappa_\Scal : \Scal \subseteq [N]\ \text{and}\ |\Scal| \geq p_1\}.
\eequationNN
The fact that this value is finite follows under the present assumptions from this section using the same arguments as in Lemma~\ref{lemma.F.psd}. Second, let $\hat\rho_1$ be defined as in Lemma~\ref{lemma.fletures.lag.result}.  Third, as in \eqref{eq.rho_lower}, but with $\kappa$ replaced by the uniform bound $\bar\kappa$ defined above and $\hat\rho_1$ included, let
\bequationNN
  \rho_\beta := \max\left\{\frac{\max\{1, \bar\kappa\}}{\sigma_{\min}}, \frac{\beta + \bar\kappa + m \bar\kappa^2}{2 \sigma_{\min}^2}, \hat\rho_1 \right\}\ \ \text{and}\ \ \rho_{\max} := \mu^{-1} \rho_\beta.
\eequationNN
Fourth, following the notation in Theorem~\ref{theo.fl.strongly.convex}, let
\begin{align*}
  \omega &:= 1 + \bar\kappa + 2\rho_{\max} \bar\kappa, \\
  \eta &:= \frac{2\sqrt{m}\bar\kappa}{\sigma_{\min}} + m \bar\kappa + 2 m \rho_{\max} \bar\kappa, \\
  M &:= M_f + \rho_{\max} M_c,
\end{align*}
and let $\bar\lambda \in \R{}_{>0}$ be any real number satisfying
\bequation\label{eq.lambdabar}
\baligned
  \bar\lambda \geq \beta + \sup_{(\rho,\Scal)} &\ \max \left\{0,\sup_{x \in \Ccal_{\Scal,R}} \lambda_1(\nabla^2 F_{\Scal,\rho}(x)) \right\} \\
  \st &\ \rho \in (\rho_\beta, \rho_{\max}] \\
  &\ \Scal \subseteq [N]\ \text{with}\ |\Scal| \geq p_1.
\ealigned
\eequation
Finally, define the values
\begin{align}
  \epsilon_* &:= \min\left\{\frac{\beta}{4\(\eta + \frac{8 M \omega}{\beta} +\(1+\frac{2}{\sigma_{\min}}\) \frac{\bar\kappa}{\sqrt{5}}\)}, \frac{R}{\omega}, \frac{\beta\bar\lambda}{4M\omega}, \frac{\bar\lambda \delta}{2\omega} \right\} \label{eq.epsstar} \\
  \epsilon_c &:= \tau_1 \xi_{\Scal_1} \label{eq.epsc} \\ \text{and}\ \ 
  c &:= \min \left\{ \frac{1}{3\kappa_1}, \frac{\alpha}{4\kappa_2}, \frac{2\beta}{9\kappa_3}, \frac{\epsilon_*}{3\tau_1}, \frac{\sqrt{5} R}{2\tau_1}, \frac{1}{\tau_2} \right\}. \label{eq.cccc}
\end{align}
It is useful to state the following technical lemma.

\blemma\label{lem.technical}
  For any $b \in \R{}_{>0}$, $\kappa \in \R{}_{>0}$, and $\rho \in \R{}_{>0}$, let $\epsilon_\rho(b,\kappa)$ denote the right-hand side of~\eqref{ineq.requirement.theorem} with $\beta$ replaced by $b$, $(\omega_\rho,\eta_\rho)$ defined as in \eqref{eq.def.eta_1,2} for the given $\kappa$, $M_\rho = M_f + \rho M_c$, and $\lambda_\rho$ replaced by $\bar\lambda$ defined in \eqref{eq.lambdabar}, i.e.,
  \bequation\label{eq.eps.rho.explicit}
    \epsilon_\rho(b,\kappa) = \min \left\{ \frac{b}{2\(\eta_\rho + \frac{4M_\rho\omega_\rho}{b} + \(1 + \tfrac{2}{\sigma_{\min}}\)\tfrac{\kappa}{\sqrt{5}}\)}, \frac{R}{\omega_\rho}, \frac{b\bar\lambda}{2M_\rho\omega_\rho}, \frac{\bar\lambda \delta}{2\omega_\rho} \right\}.
  \eequation
  Then, $\epsilon_\rho(b,\kappa)$ is nondecreasing in $b$ and nonincreasing in each of $\kappa$ and $\rho$.  Consequently, for every $\rho \in (\rho_\beta,\rho_{\max}]$, $\kappa \in (0, \bar\kappa]$, and $b \in [\beta/2,\beta]$ one has that $\epsilon_* \leq \epsilon_\rho(b,\kappa)$, while at the same time for any $\epsilon \in (0,\epsilon_\rho(b,\kappa)]$ and with
  \bequation\label{eq.trhotemp}
    t_\rho(b,\kappa) := \min \left\{\frac{1}{\bar\lambda}, \frac{1}{b - \eta_\rho \epsilon - \frac{4 M_\rho \omega_\rho}{b} \epsilon} \right\}
  \eequation
  $($see \eqref{eq.stepsize.statement}$)$ one has that $1/\bar\lambda = t_\rho(b,\kappa)$, so $[\mu/\bar\lambda,1/\bar\lambda] \subset (0,t_\rho(b,\kappa)]$.
\elemma
\bproof
  For the sake of brevity, let $C(\kappa) := \(1 + \tfrac{2}{\sigma_{\min}}\)\tfrac{\kappa}{\sqrt{5}}$. Consider first monotonicity with respect to $b$.  The first term in \eqref{eq.eps.rho.explicit} is
  \bequationNN
    \frac{b}{2\(\eta_\rho + \frac{4M_\rho\omega_\rho}{b} + C(\kappa)\)} = \frac{b^2}{2\(\(\eta_\rho + C(\kappa)\)b + 4M_\rho\omega_\rho\)},
  \eequationNN
  which is strictly increasing in $b$.  On the other hand, the second and fourth terms in \eqref{eq.eps.rho.explicit} are independent of $b$ while the third term is increasing in $b$.  Hence, overall, $\epsilon_\rho(b,\kappa)$ is nondecreasing in $b$, as claimed. Consider next monotonicity with respect to~$\kappa$.  By \eqref{eq.def.eta_1,2}, both $\omega_\rho$ and $\eta_\rho$ are increasing in $\kappa$. At the same time, $C(\kappa)$ is increasing in $\kappa$ while $M_\rho$ is independent of $\kappa$.  Consequently, each of the four terms in \eqref{eq.eps.rho.explicit} is nonincreasing in $\kappa$, and therefore so is $\epsilon_\rho(b,\kappa)$. Consider finally monotonicity with respect to~$\rho$.  Each of $\omega_\rho$, $\eta_\rho$, and $M_\rho$ is increasing in~$\rho$ while $\bar\lambda$ is independent of $\rho$, so by similar reasoning each of the four terms in \eqref{eq.eps.rho.explicit} is nonincreasing in $\rho$, and therefore so is $\epsilon_\rho(b,\kappa)$.

  Now consider arbitrary $\rho \in (\rho_\beta,\rho_{\max}]$, $\kappa \in (0,\bar\kappa]$, and $b \in [\beta/2,\beta]$.  Combining the three monotonicity properties established above, one finds that
  \bequation\label{eq.leb}
    \epsilon_\rho(b,\kappa) \geq \epsilon_{\rho_{\max}}\!\(\tfrac{1}{2}\beta,\bar\kappa\).
  \eequation
  On the other hand, since $(\omega_{\rho_{\max}},\eta_{\rho_{\max}},M_{\rho_{\max}}) = (\omega,\eta,M)$ by the definitions of $\omega$, $\eta$, and $M$, evaluating \eqref{eq.eps.rho.explicit} at $(b,\kappa,\rho) = (\beta/2,\bar\kappa,\rho_{\max})$ yields
  \bequationNN
    \epsilon_{\rho_{\max}} \(\beta/2,\bar\kappa\) = \min\left\{ \frac{\beta}{4\(\eta + \frac{8M\omega}{\beta} + \(1 + \frac{2}{\sigma_{\min}}\)\frac{\bar\kappa}{\sqrt{5}}\)}, \frac{R}{\omega}, \frac{\beta\bar\lambda}{4M\omega} \right\} = \epsilon^*,
  \eequationNN
  which with \eqref{eq.leb} shows that $\epsilon^* \leq \epsilon_\rho(b,\kappa)$, as desired.

  Lastly, for the arbitrary $\rho \in (\rho_\beta,\rho_{\max}]$, $\kappa \in (0,\bar\kappa]$, and $b \in [\beta/2,\beta]$, consider arbitrary $\epsilon \in (0,\epsilon_\rho(b,\kappa)]$. By \eqref{eq.lambdabar} one has $\bar\lambda \geq \beta$, which with $b \leq \beta$ gives
  \bequationNN
    b - \eta_\rho\epsilon - \frac{4M_\rho\omega_\rho}{b}\epsilon \leq b \leq \beta \leq \bar\lambda.
  \eequationNN
  Consequently, the first term in the minimum defining $t_\rho(b,\kappa)$ is the smaller one and $t_\rho(b,\kappa) = 1/\bar\lambda$. Thus, indeed, $1/\bar\lambda \leq t_\rho(b,\kappa)$, and since $\mu \in (0,1)$ it follows that $[\mu/\bar\lambda,1/\bar\lambda] \subseteq (0,t_\rho(b,\kappa)]$, which completes the proof.
  \qed
\eproof

We now prove our main theorem for progressive sampling.  In the theorem, for a reasonable basis of comparison, our primary measures of performance of an algorithm are the number of individual objective gradients and individual constraint Jacobians that need to be computed prior to termination with an approximate second-order stationary point of a desired user-defined accuracy.  In other words, we are concerned with the number of times that $\nabla f_i(x)$ is evaluated for some $i \in [N]$ and some $x \in \R{n}$, and the number of times that $\nabla c_i(x)$ is evaluated for some $i \in [N]$ and some $x \in \R{n}$.  One might also be interested in the numbers of individual objective function evaluations, individual constraint function evaluations, or evaluations of a least square multiplier. That said, if the methods are implemented with fixed step-size strategies, then the algorithmic complexity with respect to these quantities would be the same as for the derivative computations that we consider.

\begin{theorem}\label{theo.complexity.OLD}
  Suppose that Assumptions~\ref{ass.app.boundness}, \ref{ass.app.bounded.distribute}, \ref{ass.app.compact}, \ref{ass.app.licq}, \ref{ass.app.Lagrangian}, and \ref{ass.app.subproblems} hold.  Then, the following hold for Algorithm~\ref{alg.psm}, where it is assumed that calls to Algorithm~\ref{alg.sub.solver} for $k=1$ use either the line-search methods in \cite[Algorithms~2--3]{GoyeEfteBoum2024} or set the step sizes according to the lower bounds in \cite[Lemmas~3.2 and 3.4]{GoyeEfteBoum2024}.  $($In the following parts, we refer to $\bar\epsilon_1(\cdot,\cdot)$ and $\bar\zeta_1(\cdot,\cdot)$ defined in Lemma~\ref{lemma.fletures.lag.result}.$)$
  \benumerate
    \item[(a)] Suppose that $p_1 = N$,  $\Scal_1 = [N]$, $x_0 \in \Ccal_{[N],R}$, and $\rho_1 \geq \hat\rho_1$, where $\hat\rho_1$ is defined as in~Lemma~\ref{lemma.fletures.lag.result}. Then, for any $(\epsilon,\zeta) \in (0,\frac{\sqrt{5}}{2}R] \times (0,1]$, Algorithm~\ref{alg.psm} with initial point~$x_0$ and tolerances $(\epsilon_1,\zeta_1)$ satisfying
    \bequationNN
      \epsilon_1 \in [\mu \bar\epsilon_1(\epsilon,\zeta),\bar\epsilon_1(\epsilon,\zeta)]\ \ \text{and}\ \ \zeta_1 \in [\mu \bar\zeta_1(\epsilon,\zeta), \bar\zeta_1(\epsilon,\zeta)]
    \eequationNN
    for some $\mu \in (0,1)$ requires
    \bequation\label{eq.total.grad.eval.dir}
      \Ocal\( N \lceil \max \{\epsilon^{-2}, \zeta^{-3} \} \rceil \)
    \eequation
    individual objective gradients and individual constraint Jacobians until it terminates, yielding $x_1$ that is $(\epsilon,\zeta)$-stationary for~\eqref{prob.opt.N}.
    \item[(b)] Suppose that $p_1 < N$, $\Scal_1 \subset [N]$ with $|\Scal_1| = p_1$ satisfies $\xi_{\Scal_1} \leq c$ $($where $c$ is defined in~\eqref{eq.cccc}$)$, $x_0 \in \mathcal{C}_{\Scal_1,R}$, and the ultimate desired tolerances satisfy $(\epsilon,\zeta) \in (0,3\epsilon_c) \times \R{}_{>0}$ $($where $\epsilon_c$ is defined in~\eqref{eq.epsc}$)$.  Suppose further that the parameters of Algorithm~\ref{alg.psm} and the calls to Algorithm~\ref{alg.sub.solver} have:
    \bitemize[label=$\bullet$]
      \item $\rho_k \in (\rho_\beta,\rho_{\max}]$ for all $k \in [K]$;
      \item for each $k \in \{2,\dots,K\}$, line~4 of Algorithm~\ref{alg.sub.solver} employs a fixed step size chosen to satisfy $t_g \in [\mu/\bar\lambda,1/\bar\lambda]$;
      \item for all $k \in [K-1]$, target accuracies $(\hat\epsilon_k,\hat\zeta_k) \in \R{}_{>0} \times \R{}_{>0}$ are chosen such that $\hat\epsilon_k \in [\mu\epsilon_c,\epsilon_c]$ and $\hat\zeta_k \in [\mu\tau_2\xi_{\Scal_k},\tau_2\xi_{\Scal_k}]$;
      \item the tolerances input to Algorithm~\ref{alg.sub.solver} satisfy
      \begin{align*}
        \epsilon_1 &\in [\mu\bar\epsilon_1(\hat\epsilon_1,\hat\zeta_1), \bar\epsilon_1(\hat\epsilon_1,\hat\zeta_1)], \\
        \zeta_1 &\in [\mu\bar\zeta_1(\hat\epsilon_1,\hat\zeta_1), \bar\zeta_1(\hat\epsilon_1,\hat\zeta_1)], \\
        (\epsilon_k,\zeta_k) &= (\hat\epsilon_k/\sqrt{5},\hat\zeta_k)\ \ \text{for all}\ \ k \in \{2,\dots,K-1\}, \\ \text{and}\ \ 
        (\epsilon_K,\zeta_K) &= (\epsilon/\sqrt{5},\zeta).
      \end{align*}
    \eitemize
    Then, with
    \bequation\label{eq.B.Tbar}
      B := \frac{4}{4 - \frac{\mu\beta}{2\bar\lambda}} > 1,\ \ \overline{T} := \left\lceil \log_B \(\frac{3\sqrt{5}\omega}{\mu}\) \right\rceil,
    \eequation
    $(u_{\Scal_1,1},u_{\Scal_1,2})$ defined as in Lemma~\ref{lemma.fletures.lag.result}, and
    \bequation\label{eq.C1}
      C_1 := \max\left\{ \frac{u_{\Scal_1,1}}{\mu^{4}\xi_{\Scal_1}^2} \max\left\{\frac{5}{\tau_1^2},\frac{16\tau_{\Scal_1}^2}{\tau_2^2}\right\}, \frac{8u_{\Scal_1,2}}{\mu^{6}(\tau_2\xi_{\Scal_1})^{3}} \right\},
    \eequation
    Algorithm~\ref{alg.psm} terminates with $x_K$ satisfying
    \begin{align*}
      \|\nabla L(x_K,y(x_K))\| &\leq \epsilon \\ \text{and}\ \ 
      d^T\nabla^2_{xx}L(x_K,y(x_K))d &\geq \tfrac{1}{2} \beta \|d\|^2\ \ \text{for all}\ \ d \in \Null(\nabla c(x_K)^T),
    \end{align*}
    so that $x_K$ is $(\epsilon,\zeta')$-stationary for~\eqref{prob.opt.N} for every $\zeta' \in \R{}_{\geq 0}$.  Moreover, the total number of individual objective gradients and individual constraint Jacobians computed before termination is at most
    \bequation\label{eq.total}
      p_1 C_1 + \overline{T} \frac{\theta(N - p_1)}{\theta - 1} + N \left\lceil \log_B \(\frac{3\sqrt{5} \omega \epsilon_c}{\epsilon} \)\right\rceil,
    \eequation
    where the first two terms are independent of $\epsilon$ and $\zeta$.
  \eenumerate
\end{theorem}
\bproof
  Consider part (a). By Lemma~\ref{lemma.fletures.lag.result}, there exists a pair $(u_{[N],1},u_{[N],2}) \in \R{}_{>0} \times \R{}_{>0}$ such that, for any $(\epsilon,\zeta) \in (0,\frac{\sqrt{5}}{2}R] \times (0,1]$ and tolerances $(\epsilon_1,\zeta_1)$ satisfying the conditions of part (a), Algorithm~\ref{alg.psm} requires at most
  \bequationNN
    T_{\Scal_1} = \max\{u_{[N],1} \epsilon_1^{-2}, u_{[N],2} \zeta_1^{-3}\}
  \eequationNN
  iterations until $x_1$ is produced, which is $(\epsilon,\zeta)$-stationary, as desired.  Further, since $\epsilon_1 \geq \mu \bar\epsilon_1(\epsilon,\zeta) = \mu \min\{\frac{\sqrt{5}\epsilon}{5}, \frac{\zeta}{4\tau_{[N]}}\}$ and $\zeta_1 \geq \mu \bar\zeta_1(\epsilon,\zeta) \geq \mu \frac{\zeta}{2}$, one finds
  \begin{align*}
    T_{\Scal_1}
      &\leq \max\{5 \mu^{-2} u_{[N],1} \epsilon^{-2}, 16\tau_{[N]}^2 \mu^{-2} u_{[N],1} \zeta^{-2}, 8 \mu^{-3} u_{[N],2} \zeta^{-3}\} \\
      &\leq \max\{5 \mu^{-2} u_{[N],1} \epsilon^{-2}, 16\tau_{[N]}^2 \mu^{-2} u_{[N],1} \zeta^{-3}, 8 \mu^{-3} u_{[N],2} \zeta^{-3}\}.
  \end{align*}
  Thus, the claim follows since each iteration the call to Algorithm~\ref{alg.sub.solver} requires $N$ individual objective gradients and $N$ individual constraint Jacobians.

  Now consider part (b).  Throughout the proof of this part, for each $k \in [K]$ let $\kappa_{\Scal_k} \leq \bar\kappa$ denote a value satisfying \eqref{ineq.bound.bar.kappa} with respect to \eqref{prob.opt.S} and define the pair $(\alpha_{\Scal_k},\beta_{\Scal_k})$ in the manner of Lemma~\ref{lem.strongMorseextend}.  Since $p_{k+1} = \min\{\theta p_k, N\}$ for all $k \in [K]$ and $K = \lceil \log_\theta (N/p_1)\rceil + 1$, one finds that $p_k = \theta^{k-1}p_1 < N$ for all $k \in [K-1]$ and $p_K = N$.  Consequently, $\xi_{\Scal_k} > 0$ for all $k \in [K-1]$ while $\xi_{\Scal_K} = 0$ so $\alpha_{\Scal_K} = \alpha$ and $\beta_{\Scal_K} = \beta$. In addition, since $\xi_{\Scal_1} \leq c$ and, by \eqref{eq.key_value} and $|\Scal_k| \geq p_1$, one has $\xi_{\Scal_k} \leq \xi_{\Scal_1}$ for all $k \in [K]$, it follows from the definition of $c$ that \eqref{eq.theorem2.S} holds.  Hence, by Lemma~\ref{lem.strongMorseextend}, $\beta_{\Scal_k} \in [\beta/2,\beta]$ for all $k \in [K]$. Finally, observe from \eqref{eq.epsc} and \eqref{eq.cccc} and the conditions of part (b) that
  \bequation\label{eq.basin}
    3\epsilon_c = 3\tau_1\xi_{\Scal_1} \leq 3\tau_1 c \leq \epsilon^* \leq \frac{R}{\omega} \leq R.
  \eequation

  We now claim that, for any $k \in \{2,\dots,K\}$ and any $\epsilon' \in (0,3\epsilon_c)$, if
  \bequation\label{eq.invariant}
    \begin{aligned}
      \|\nabla L_{\Scal_k}(x_{k-1},y_{\Scal_k}(x_{k-1}))\| &\leq 3\epsilon_c \\
      \text{and}\ \ d^T\nabla^2_{xx}L_{\Scal_k}(x_{k-1},y_{\Scal_k}(x_{k-1}))d &\geq \beta_{\Scal_k}\|d\|^2 \\ 
      \text{for all}\ \ d &\in \Null(\nabla c_{\Scal_k}(x_{k-1})^T),
    \end{aligned}
  \eequation
  then Algorithm~\ref{alg.sub.solver}, called with inputs $x_{k-1}$, $\Scal_k$, $\rho_k$, and $(\epsilon_k,\zeta_k)$ replaced by tolerances $(\epsilon'/\sqrt{5},\zeta_k)$ and with the fixed step size $t_g$, reduces to gradient descent applied to minimize $F_{\Scal_k,\rho_k}$, terminates in at most $\lceil \log_B(3\sqrt{5}\omega/\epsilon')\epsilon_c\rceil$---more precisely, at most $\lceil \log_B(3\sqrt{5}\omega\epsilon_c/\epsilon')\rceil$---iterations, and returns a point $x_k$ satisfying \eqref{eq.epsbeta} with respect to \eqref{prob.opt.S} with $(\epsilon,\beta)$ replaced by $(\epsilon',\beta_{\Scal_k}/2)$.

  To prove this claim, we verify the conditions of Theorem~\ref{theo.fl.strongly.convex} for~\eqref{prob.opt.S} with $\rho = \rho_k$, $\beta = \beta_{\Scal_k}$, $\kappa = \kappa_{\Scal_k}$, $\epsilon = 3\epsilon_c$, $\lambda_\rho = \bar\lambda$, and $t = t_g$. First, the threshold in \eqref{eq.rho_lower} is increasing in each of $\beta$ and $\kappa$, so with $\beta_{\Scal_k} \leq \beta$ and $\kappa_{\Scal_k} \leq \bar\kappa$ its value is at most $\rho_\beta$, meaning that \eqref{eq.rho_lower} holds since $\rho_k > \rho_\beta$.  Second, by \eqref{eq.lambdabar} one has $\bar\lambda \geq \beta + \sup_{x \in \mathcal{C}_{\Scal_k,R}} \lambda_1(\nabla^2F_{\Scal_k,\rho_k}(x)) \geq \beta_{\Scal_k} + \lambda_{R,\rho_k}$, so $\bar\lambda$ is valid in place of $\lambda_\rho$.  Third, by Lemma~\ref{lem.technical} with $(b,\kappa,\rho) = (\beta_{\Scal_k},\kappa_{\Scal_k},\rho_k)$---for which $b \in [\beta/2,\beta]$, $\kappa \leq \bar\kappa$, and $\rho \in (\rho_\beta,\rho_{\max}]$ hold---and by \eqref{eq.basin}, one finds $3\epsilon_c \leq \epsilon^* \leq \epsilon_{\rho_k}(\beta_{\Scal_k},\kappa_{\Scal_k})$, so \eqref{ineq.requirement.theorem} holds. In addition, $3\epsilon_c \in (0,R)$ by \eqref{eq.basin}, as required in Lemma~\ref{lemma.F.psd}. Fourth, by Lemma~\ref{lem.technical} one has $t_g \leq 1/\bar\lambda \leq t_{\rho_k}$, so \eqref{eq.stepsize.statement} holds.  Finally, \eqref{eq.invariant} states that $x_{k-1}$ satisfies \eqref{eq.epsbeta} with $\epsilon = 3\epsilon_c$ and $\beta = \beta_{\Scal_k}$, and $\epsilon' \in (0,3\epsilon_c)$. Consequently, Theorem~\ref{theo.fl.strongly.convex} applies.  By the parameter correspondence \eqref{goyens.alg.parameter}, the tolerance $\epsilon'/\sqrt{5}$ and step size $t_g$ input to Algorithm~\ref{alg.sub.solver} are exactly those in \eqref{goyens.alg.parameter}, so Algorithm~\ref{alg.sub.solver} reduces to the gradient-descent iteration \eqref{eq.fl.gd}; indeed, by \eqref{eq.beta.eps.rho} and \eqref{ineq.strong.convex.2} one has $\lambda_n(\nabla^2F_{\Scal_k,\rho_k}(x^i_{k-1})) \geq \beta_{\Scal_k}/2 > 0$ for all $i$, so the condition in line~5 of Algorithm~\ref{alg.sub.solver} never holds.  Moreover, by \eqref{eq.T}, the number of iterations performed is at most
  \bequationNN
    \left\lceil \log_{\frac{4}{4-\beta_{\Scal_k}t_g}} \left(\frac{3\sqrt{5} \omega_{\rho_k}\epsilon_c}{\epsilon'}\right)\right\rceil \leq \left\lceil \log_B\left(\frac{3\sqrt{5} \omega\epsilon_c}{\epsilon'}\right)\right\rceil,
  \eequationNN
  where the inequality follows since $\omega_{\rho_k} \leq \omega$ (as $\rho_k \leq \rho_{\max}$ and $\kappa_{\Scal_k} \leq \bar\kappa$), since $3\sqrt{5}\omega\epsilon_c/\epsilon' > 1$, and since the function $4/(4-a)$ is increasing over $a \in [0,4)$, which with $\beta_{\Scal_k}t_g \geq \beta \mu/(2\bar\lambda)$ shows that $4/(4-\beta_{\Scal_k}t_g) \geq B$.  Lastly, by \eqref{eq.beta'}, the returned point $x_k$ satisfies \eqref{eq.epsbeta} with $(\epsilon,\beta)$ replaced by $(\epsilon',\beta')$ where $\beta' \geq \beta_{\Scal_k}/2 \geq \beta/4 > 0$, which proves the claim.

  Now consider $k=1$. By~\eqref{eq.cccc} one has $\hat\epsilon_1 \leq \epsilon_c = \tau_1\xi_{\Scal_1} \leq \tau_1 c \leq \tfrac{\sqrt{5}}{2}R$ and $\hat\zeta_1 \leq \tau_2\xi_{\Scal_1} \leq \tau_2 c \leq 1$, so $(\hat\epsilon_1,\hat\zeta_1) \in (0,\tfrac{\sqrt{5}}{2}R] \times (0,1]$ satisfies the conditions of $(\bar\epsilon,\bar\zeta)$ in Lemma~\ref{lemma.fletures.lag.result}.  In addition, $x_0 \in \mathcal{C}_{\Scal_1,R}$ and $\rho_1 > \rho_\beta \geq \hat\rho_1$.  Hence, by Lemma~\ref{lemma.fletures.lag.result}, the call to Algorithm~\ref{alg.sub.solver} with inputs $x_0$, $\Scal_1$, $\rho_1$, and $(\epsilon_1,\zeta_1)$ terminates in at most $T_{\Scal_1} = \max\{u_{\Scal_1,1}\epsilon_1^{-2},u_{\Scal_1,2}\zeta_1^{-3}\}$ iterations with $x_1$ being $(\hat\epsilon_1,\hat\zeta_1)$-stationary with respect to \eqref{prob.opt.S} with $k = 1$.  Furthermore,
  \begin{align*}
    \epsilon_1 \geq \mu\bar\epsilon_1(\hat\epsilon_1,\hat\zeta_1) &= \mu\min\{\sqrt{5}\hat\epsilon_1/5,\hat\zeta_1/(4\tau_{\Scal_1})\} \\
    &\geq \mu^2\xi_{\Scal_1}\min\{\sqrt{5}\tau_1/5,\tau_2/(4\tau_{\Scal_1})\} \\
    \text{and}\ \ \zeta_1 \geq \mu\bar\zeta_1(\hat\epsilon_1,\hat\zeta_1) &\geq \mu\hat\zeta_1/2 \geq \mu^2\tau_2\xi_{\Scal_1}/2,
  \end{align*}
  so one finds from \eqref{eq.C1} that $T_{\Scal_1} \leq C_1$.  Since each iteration of Algorithm~\ref{alg.sub.solver} in this case requires $p_1$ individual objective gradients and $p_1$ individual constraint Jacobians, $k = 1$ requires at most $p_1C_1$ of each, the first term in \eqref{eq.total}.

  We now prove by induction that, for each $k \in \{2,\dots,K\}$, the point $x_{k-1}$ satisfies \eqref{eq.invariant} and, if $k \leq K-1$, that Algorithm~\ref{alg.sub.solver} returns $x_k$ that is $(\hat\epsilon_k,\hat\zeta_k)$-stationary with respect to \eqref{prob.opt.S} in at most $\bar{T}$ iterations. For the base case, recall from the previous paragraph that $x_1$ is $(\hat\epsilon_1,\hat\zeta_1)$-stationary with respect to \eqref{prob.opt.S} with $k = 1$, where $\hat\epsilon_1 \leq \epsilon_c = \tau_1\xi_{\Scal_1}$ and $\hat\zeta_1 \leq \tau_2\xi_{\Scal_1}$, so that \eqref{eq.tolorence} holds for $k = 1$.  Since \eqref{eq.theorem2.S} holds by previous arguments, Lemma~\ref{lem.start_for_next} with $k = 1$ yields \eqref{ineq.init.next}, which is exactly \eqref{eq.invariant} for $k = 2$ since $3\tau_1\xi_{\Scal_1} = 3\epsilon_c$. Next, for the inductive step, consider $k \in \{2,\dots,K-1\}$ and suppose that $x_{k-1}$ satisfies \eqref{eq.invariant}.  Since $\epsilon_k = \hat\epsilon_k/\sqrt{5}$ with $\hat\epsilon_k \in [\mu\epsilon_c,\epsilon_c]$, one has $\hat\epsilon_k \in (0,3\epsilon_c)$, so Theorem~\ref{theo.fl.strongly.convex} applies with $\epsilon' = \hat\epsilon_k$.  Hence Algorithm~\ref{alg.sub.solver} returns $x_k$ in at most
  \bequationNN
    \left\lceil \log_B\left(\frac{3\sqrt{5}\,\omega\epsilon_c}{\hat\epsilon_k}\right)\right\rceil \leq \left\lceil \log_B\left(\frac{3\sqrt{5}\,\omega}{\mu}\right)\right\rceil = \bar{T}
  \eequationNN
  iterations, and $x_k$ satisfies \eqref{eq.epsbeta} with respect to \eqref{prob.opt.S} with $(\epsilon,\beta)$ replaced by $(\hat\epsilon_k,\tfrac{1}{2}\beta_{\Scal_k})$.  In particular, $\|\nabla L_{\Scal_k}(x_k,y_{\Scal_k}(x_k))\| \leq \hat\epsilon_k$ and, since one trivially finds $\beta_{\Scal_k}/2 > 0 \geq -\hat\zeta_k$, it follows that $x_k$ is $(\hat\epsilon_k,\hat\zeta_k)$-stationary with respect to \eqref{prob.opt.S}.  Now since $\hat\epsilon_k \leq \epsilon_c = \tau_1\xi_{\Scal_1}$ and $\hat\zeta_k \leq \tau_2\xi_{\Scal_k}$, the tolerances satisfy \eqref{eq.tolorence}, so Lemma~\ref{lem.start_for_next} yields \eqref{ineq.init.next}, which is \eqref{eq.invariant} with $k$ replaced by $k+1$.  This completes the induction.  Overall, since each iteration of Algorithm~\ref{alg.sub.solver} in iteration $k$ requires $p_k$ individual objective gradients and $p_k$ individual constraint Jacobians, and since $p_k = \theta^{k-1}p_1$ for all $k \in [K-1]$, the total of these evaluations required over iterations $k \in \{2,\dots,K-1\}$ is at most
  \bequationNN
    \bar{T}\sum_{k=2}^{K-1} p_k = \bar{T}\,p_1\sum_{j=1}^{K-2}\theta^j = \bar{T}\,\frac{\theta(p_{K-1} - p_1)}{\theta - 1} \leq \bar{T}\,\frac{\theta(N - p_1)}{\theta - 1},
  \eequationNN
  which is the second term in \eqref{eq.total}.  (We remark specifically that if $K = 2$, then this sum is empty and the bound holds trivially.)

  Finally, consider $k = K$.  By the previous induction, $x_{K-1}$ satisfies \eqref{eq.invariant} with $k = K$, where $\beta_{\Scal_K} = \beta$ and $L_{\Scal_K} = L$ by previous arguments. Since $\epsilon_K = \epsilon/\sqrt{5}$ with $\epsilon \in (0,3\epsilon_c)$, Theorem~\ref{theo.fl.strongly.convex} applies with $\epsilon' = \epsilon$, so Algorithm~\ref{alg.sub.solver} returns $x_K$ in at most $\lceil \log_B(3\sqrt{5}\omega\epsilon_c/\epsilon)\rceil$ iterations, and $x_K$ satisfies \eqref{eq.epsbeta} with respect to \eqref{prob.opt.N} with $(\epsilon,\beta)$ replaced by $(\epsilon,\beta/2)$, i.e., $\|\nabla L(x_K,y(x_K))\| \leq \epsilon$ and $d^T\nabla^2_{xx}L(x_K,y(x_K))d \geq \tfrac{1}{2}\beta\|d\|^2$ for all $d \in \Null(\nabla c(x_K)^T)$.  Since $\tfrac{1}{2}\beta > 0 \geq -\zeta'$ for any $\zeta' \in \R{}_{\geq 0}$, it follows that $x_K$ is $(\epsilon,\zeta')$-stationary for \eqref{prob.opt.N} for every $\zeta' \in \R{}_{\geq 0}$, as claimed.  Since each iteration in this case requires $N$ individual objective gradients and $N$ individual constraint Jacobians, iteration $k = K$ requires at most $N\lceil \log_B(3\sqrt{5}\omega\epsilon_c/\epsilon)\rceil$ of each, which is the third term in \eqref{eq.total}.  Summing the bounds from the previous arguments completes the proof.
  \qed
\eproof

\section{Conclusion}\label{sec.conclusion}

We have proved that a recently proposed Gradient-Eigenstep algorithm for solving equality constrained continuous optimization problems, which is based on minimizing Fletcher's augmented Lagrangian function, reduces to a local-linearly convergent gradient-descent method in the vicinity of certain approximate second-order stationary points when the step size is sufficiently small and the penalty parameter is sufficiently large. We have also leveraged this result to show that a progressive sampling strategy for solving sample-average-type problems can obtain an improved worst-case sample complexity bound as compared to an approach that solves a full-sample problem directly.

\section*{Data Availability}

We do not analyze or generate any datasets because our work proceeds within a theoretical and mathematical approach.

\bibliographystyle{plain}
\bibliography{references}

@article{GoyeEfteBoum2024,
  title = {Computing {{Second-Order Points Under Equality Constraints}}: {{Revisiting Fletcher}}'s {{Augmented Lagrangian}}},
  author = {Goyens, Florentin and Eftekhari, Armin and Boumal, Nicolas},
  year = {2024},
  journal = {Journal of Optimization Theory and Applications},
  volume = {201},
  number = {3},
  pages = {1198--1228}
}

@book{Stew1973,
  title={Introduction to Matrix Computations},
  author={Stewart, G. W.},
  series={Computer Science and Applied Mathematics},
  year={1973},
  publisher={Academic Press}
}

@article{Stew1977,
 author = {G. W. Stewart},
 journal = {SIAM Review},
 number = {4},
 pages = {634--662},
 publisher = {Society for Industrial and Applied Mathematics},
 title = {On the Perturbation of Pseudo-Inverses, Projections and Linear Least Squares Problems},
 volume = {19},
 year = {1977}
}

@article{Stew1979,
title = {A note on the perturbation of singular values},
journal = {Linear Algebra and its Applications},
volume = {28},
pages = {213-216},
year = {1979},
author = {G. W. Stewart}
}

@article{gallier2019linear,
  title={Linear algebra for computer vision, robotics, and machine learning},
  author={Gallier, Jean and Quaintance, Jocelyn},
  journal={University of Pennsylvania},
  year={2019}
}

@article{fujiwara1982morse,
  title={Morse programs: a topological approach to smooth constrained optimization},
  author={Fujiwara, Okitsugu},
  journal={Mathematics of Operations Research},
  volume={7},
  number={4},
  pages={602--616},
  year={1982},
  publisher={INFORMS}
}

@book{milnor1963morse,
  title={Morse theory},
  author={Milnor, John Willard},
  number={51},
  year={1963},
  publisher={Princeton university press}
}

@book{guillemin1975topology,
isbn = {0132126052},
language = {eng},
lccn = {74004115},
publisher = {Prentice-Hall},
title = {Differential topology },
year = {1974},
author = {Guillemin, V. and Pollack, Alan.},
address = {Englewood Cliffs, N.J},
booktitle = {Differential topology},
}

@unpublished{CurtGuoRobi2025,
author = {Frank E. Curtis and Lingjun Guo and Daniel P. Robinson},
title = {Progressively Sampled Equality-Constrained Optimization},
note = {arXiv 2510.00417},
year = {2025}
}

@article{estrin2020,
  title={Implementing a smooth exact penalty function for equality-constrained nonlinear optimization},
  author={Estrin, Ron and Friedlander, Michael P and Orban, Dominique and Saunders, Michael A},
  journal={SIAM Journal on Scientific Computing},
  volume={42},
  number={3},
  pages={A1809--A1835},
  year={2020},
  publisher={SIAM}
}

@misc{fletcher1970class,
  title={A class of methods for nonlinear programming with termination and convergence properties},
  author={Fletcher, Roger},
  journal={Integer and nonlinear programming},
  volume={157175},
  year={1970},
  publisher={North-Holland Amsterdam}
}

@incollection{DiPillo1994,
  author    = {Di Pillo, Gianni},
  title     = {Exact Penalty Methods},
  booktitle = {Algorithms for Continuous Optimization},
  pages     = {209--253},
  year      = {1994},
  publisher = {Springer},
  address   = {Dordrecht},
  doi       = {10.1007/978-94-009-0369-2_8}
}

@article{DiPilloGrippo1986,
  author    = {Di Pillo, Gianni and Grippo, Luigi},
  title     = {An Exact Penalty Function Method with Global Convergence Properties for Nonlinear Programming Problems},
  journal   = {Mathematical Programming},
  volume    = {36},
  number    = {1},
  pages     = {1--18},
  year      = {1986},
  doi       = {10.1007/BF02591986}
}

@article{DiPilloGrippo1989,
  author    = {Di Pillo, Gianni and Grippo, Luigi},
  title     = {Exact Penalty Functions in Constrained Optimization},
  journal   = {SIAM Journal on Control and Optimization},
  volume    = {27},
  number    = {6},
  pages     = {1333--1360},
  year      = {1989},
  doi       = {10.1137/0327068}
}

@book{nocedal2006numerical,
  title={Numerical optimization},
  author={Nocedal, Jorge and Wright, Stephen J},
  year={2006},
  publisher={Springer}
}

\end{document}